\documentstyle[12pt,draft]{article}

 \newcounter{abceqn}

 \newcounter{abcfig}

\newcommand{\LL}{{\mathcal L}}

\newcommand{\be}{\beta}
\newcommand{\na}{\nabla}
\newcommand{\om}{\omega}
\newcommand{\Ga}{\Gamma}
\newcommand{\Om}{\Omega}

\newcommand{\Z}{Z^2/\{0\}}
\newcommand{\pa}{\partial}

\newcommand{\tH}{\tilde{H}}

\newcommand{\e}{\epsilon}

\newcommand{\k}{\kappa}
\newcommand{\ga}{\gamma}

\newcommand{\dl}{\delta}
\newcommand{\Dl}{\Delta}

\newcommand{\ra}{\rightarrow}
\newcommand{\al}{\alpha}
\newcommand{\sg}{\sigma}
\newcommand{\Sg}{\Sigma}

\newcommand{\la}{\lambda}
\newcommand{\tla}{\tilde{\lambda}}

\newcommand{\nid}{\noindent}

\newcommand{\hk}{\hat{k}}

\newcommand{\C}{{\cal C}}

\renewcommand{\theequation}{\thesection.\arabic{equation}}
\newcommand{\eqnsection}[1]{
	\section{#1}
	\setcounter{equation}{0}
	\renewcommand{\theequation}{\thesection.\arabic{equation}}
	\setcounter{figure}{0}
	\renewcommand{\thefigure}{\thesection.\arabic{figure}}
	\setcounter{remark}{0}
	\renewcommand{\theremark}{\thesection.\arabic{remark}}
	\setcounter{theorem}{0}
	\renewcommand{\thetheorem}{\thesection.\arabic{theorem}}
	\setcounter{lemma}{0}
	\renewcommand{\thelemma}{\thesection.\arabic{lemma}}
}

\title{{\bf On 2D Euler Equations: \Large Part II. Lax Pairs and 
Homoclinic Structures}}

\author{ \\ \\ \\ \\ 
Yanguang (Charles)\ \ Li  \thanks{This work is supported 
by the AMS Centennial Fellowship, and the Guggenheim Fellowship.}
\\  \\  \\ Department of Mathematics 
 \\ \\ University of Missouri \\ \\ 
Columbia, MO 65211 \\ \\ e-mail: cli@math.missouri.edu \\ \\ 
phone: 573-884-0622 \\ \\ fax: 573-882-1869}

\date{\today}

\renewcommand{\theequation}{\thesection.\arabic{equation}}

\begin{document}
\bibliographystyle{plain}
\maketitle
\newpage
\begin{abstract}    
In Part I \cite{Li00} of our study on 2D Euler equation, 
we established the spectral theorem for a linearized 2D Euler equation. 
We also computed the point spectrum through continued fractions, and 
identified the eigenvalues with nonzero real parts.

In this Part II of our study, first we discuss the Lax pairs for both
2D and 3D Euler equations. The existence of Lax pairs suggests that the 
hyperbolic foliations of 2D and 3D Euler equations may be degenerate, 
i.e., there exist homoclinic structures. Then we investigate the question 
on the degeneracy v.s. nondegeneracy of the hyperbolic foliations for 
Galerkin truncations of 2D Euler equation. In particular, for a Galerkin 
truncation, we have computed 
the explicit representation of the hyperbolic foliation which is of the 
degenerate case, i.e., figure-eight case. We also study the robustness 
of this degeneracy for a so-called dashed-line model
through higher order Melnikov functions. The first 
order and second order Melnikov functions are all identically zero,
which indicates that the degeneracy is relatively robust. The study in 
this paper serves a clue in searching for homoclinic structures for 2D 
Euler equation. The recent breakthrough result \cite{Li00c} of mine 
on the existence of a Lax pair for 2D Euler equation, strongly supports 
the possible existence of homoclinic structures for 2D Euler equation. 


MSC: 76 34  35  37 \ .

\end{abstract}

\newtheorem{lemma}{Lemma}
\newtheorem{theorem}{Theorem}
\newtheorem{corollary}{Corollary}
\newtheorem{remark}{Remark}
\newtheorem{definition}{Definition}
\newtheorem{proposition}{Proposition}
\newtheorem{assumption}{Assumption}

\newpage
\tableofcontents



\newpage
\renewcommand{\AA}{{\mathcal A}}
\newcommand{\BB}{{\mathcal B}}

\newcommand{\DD}{{\mathcal D}}
\newcommand{\EE}{{\mathcal E}}
\newcommand{\FF}{{\mathcal F}}

\newcommand{\arc}{\mathop{\rm arc}\nolimits}
\newcommand{\diag}{\mathop{\rm diag}\nolimits}

\newcommand{\sech}{\mathop{\rm sech}\nolimits}
\newcommand{\sign}{\mathop{\rm sign}\nolimits}

\newcommand{\BD}{{B\"acklund-Darboux transformations}}

\eqnsection{Introduction}

In a series of works (\cite{Li00}, the current article, and \cite{LLS00}), 
we plan to build a dynamical system theory on 2D turbulence. The governing 
equation that we are interested in is the incompressible 2D Navier-Stokes 
equation under periodic boundary conditions. We are particularly 
interested in investigating the dynamics of 2D Navier-Stokes equation in 
the infinite Reynolds number limit and of 2D Euler equation. Our approach 
is different from many other studies on 2D Navier-Stokes equation in which 
one starts with Stokes equation to prove results on 2D Navier-Stokes 
equation for small Reynolds number. In our studies, we start with 2D Euler 
equation and view 2D Navier-Stokes equation for large Reynolds number as 
a (singular) perturbation of 2D Euler equation. 2D Euler equation is a 
Hamiltonian system with infinitely many Casimirs. To understand the nature 
of turbulence, we start with investigating the hyperbolic structure of 2D 
Euler equation. We are especially interested in investigating the 
degeneracy v.s. nondegeneracy nature of the hyperbolic foliation. Degeneracy 
means the coincidence of center-unstable and center-stable manifolds
or unstable and stable foliations, i.e.,
figure-eight structure. If all the hyperbolic foliations of 2D Euler 
equation are of figure-eight, then there is no turbulent dynamics for  
2D Euler equation. 

The recent breakthrough result \cite{Li00c} of mine on the existence 
of a Lax pair for 2D Euler equation makes my idea in the last paragraph 
even more realistic. The philosophical significance of the existence of a 
Lax pair for 2D Euler equation by the author, and more recently for 3D 
Euler equation by Steve Childress \cite{Chi00}, is beyond the particular 
project undertaken here. If one defines integrability of an equation 
by the existence of a Lax pair, then both 2D and 3D Euler equations 
are integrable. More importantly, both 2D and 3D Navier-Stokes equations 
at high Reynolds numbers are near integrable systems. Such a point of view 
changes our old ideology on Euler and Navier-Stokes equations.

In \cite{Li00}, we studied a linearized 2D Euler equation at a fixed point. 
The linear system decouples into infinitely many one-dimensional invariant 
subsystems. The essential spectrum of each invariant subsystem is a band 
of continuous spectrum on the imaginary axis. Only finitely many these 
invariant subsystems have point spectra. The point spectra can be computed 
through continued fractions. Examples show that there are indeed 
eigenvalues with positive and negative real parts. Thus, there is linear 
hyperbolicity. To understand the hyperbolic foliations of 2D Euler 
equation, especially the nature of degeneracy v.s. nondegeneracy of the 
hyperbolic foliations, we first study a Galerkin truncation. The Galerkin 
truncation presented in this article allows us to compute the hyperbolic 
foliation explicitly. We call this Galerkin truncation the 
{\em{dashed-line model}}. The dashed-line model contains a control 
parameter $\e \in [0,1]$. When the control parameter vanishes, the model 
decouples into 
a sequence of five-dimensional invariant subsystems. The eigenvalues 
of one of such invariant subsystems are ``good'' approximations of 
the eigenvalues of the original linearized 2D 
Euler equation. For this five-dimensional invariant subsystem, the 
hyperbolic foliation can be calculated explicitly with elegant 
formula representation. The hyperbolic foliation has the picture 
form of a ``lip''. Both the upper and the lower halves of the ``lip'' 
are two-dimensional ellipsoidal surfaces. 
Orbits in each half of the ``lip'' are heteroclinic orbits spiraling 
into two fixed points as time approaches positive and negative infinities,
and having a ``turning point'' (``perversion'' \cite{GT98a} \cite{GT96} 
\cite{GT97a} \cite{GT97b} \cite{GT97c} \cite{GT98b} \cite{GT98c}) in the 
middle. We call these orbits {\em{perverted heteroclinic orbits}}. 
In terms of physical variables, such heteroclinic orbits represent 
heteroclinic profile-evolutions approaching two different stationary 
profiles as time approaches positive and negative infinities. Recently, 
there are more and more interests in studying heteroclinic and homoclinic 
orbits in fluids \cite{CW95} \cite{Cra96} \cite{CG97}.
For the two-dimensional Kelvin-Helmholtz problem, Craig 
and Groves found a homoclinic orbit approaching a temporally periodic 
profile as time approaches positive and negative infinities for a 
fourth order normal form \cite{CG97}. 
Numerical calculation indicates that the eigenvalues with zero real parts 
of the dashed-line model when the control parameter vanishes, form four 
narrow bands. When the control parameter is increased from $0$ to $1$, 
the sizes of these four bands increase and congregate into the one band 
continuous spectrum of the original linearized 2D Euler equation. In 
this sense, the 
eigenvalues with zero real parts of the dashed-line model  
approximate the continuous spectrum of the original linearized 2D Euler 
equation. We are interested in investigating the degeneracy v.s. 
nondegeneracy nature of the hyperbolic foliation for 2D Euler equation. 
As a first step, we will study the persistence of the ``lip'' structure 
for the dashed-line model through Melnikov functions. It turns out that 
both the first order and the second order Melnikov functions are 
identically zero. 

There have been a lot of works on Melnikov functions since the
original work of V. K. Melnikov \cite{Mel63}.  It has been proved
rigorously that a Melnikov function is the leading order term of
a certain signed distance between center-unstable and
center-stable manifolds, see for example \cite{GH83} \cite{Wig88}
\cite{CY92} \cite{LM97}. Typically, a Melnikov function is 
a temporal integral from negative infinity to positive infinity  of the 
inner product of the gradient of certain invariant with the perturbation 
term evaluated along an unperturbed heteroclinic or homoclinic orbit, 
and this integrand often has exponential decay property as time approaches 
positive and negative infinities. The importance of Melnikov functions is 
to determine the intersection between the center-unstable and the 
center-stable manifolds. The argument is as follows. Often the Melnikov 
function is easily computable. First one seeks zeros for the Melnikov 
functions; then implicit function theorem implies nearby zeros for the 
signed distance. The intersection between  the center-unstable and the 
center-stable manifolds implies the existence of orbits homoclinic to 
the center manifold.
Other variety of Melnikov functions include subharmonic Melnikov 
functions \cite{GH83}, which are temporal integrals over the finite 
period intervals, and exponentially small Melnikov functions \cite{Gel96}, etc.
In some special
circumstance, the Melnikov function can be identically zero as a 
function of parameters.  In such case, one needs to calculate the 
next leading order term of the signed distance, called the {\em{second order 
Melnikov function}}. If the ($n-1$)-th order Melnikov function is 
identically zero as a function of the parameters, the {\em{n-th order 
Melnikov function}} serves as the leading order term of the signed 
distance. There have been some works on
second order Melnikov functions, see for example \cite{Rou86} \cite{Rou89} 
\cite{Yua94} \cite{Rom95} \cite{Dan96} \cite{IP99}.  
In this paper, we are going to derive high order Melnikov functions without
rigorous justification for the dashed-line model, and calculate them to 
study the persistence of the figure-eight structure.

Our purpose of introducing the dashed-line model is for investigating the 
hyperbolic foliation of 2D Euler equation, which is different from those 
of other models. For example, the 
shell model \cite{BLLP95} \cite{SKL95}
models the interaction of energy in different shells in the
spectral space. The difference between the dashed-line model 
and the shell model can be summarized as follows: 1. The nonlinear 
terms in the dashed-line model are derived from Galerkin truncations, while 
the nonlinear terms in the shell model are postulated. On the other hand, 
both the dashed-line model and the shell model conserve the energy and the 
enstrophy. 2. The dashed-line model captures the eigenvalues of the linearized 
2D Euler equation, while the shell model has no such information at all.
3. The dashed-line model is claimed to model the dynamics of the 2D 
Euler equation in a neighborhood of a line of fixed points and hopefully 
also the global dynamics. The variables in the dashed-line model have 
precise physical meanings. The shell model is claimed to model the 
transfer of energy in different spectral shells. The variables in the shell 
model are artificial.

The article is organized as follows: In Section 2, we discuss the Lax pairs 
for both 2D and 3D Euler equations. In Section 3, we will introduce the 
dashed-line model and discuss its rationality. In Section 4, we will compute 
explicitly the figure-eight structure for the dashed-line model at the 
special parameter value $\e =0$. In Section 5, we study the persistence 
of the figure-eight structure when $\e \neq 0$ through higher order 
Melnikov functions. Section 6 is the conclusion.

\eqnsection{Lax Pairs for 2D and 3D Euler Equations}

In this section, we are going to discuss the most recent results on the 
Lax pairs of 2D and 3D Euler equations.

\subsection{The Lax Pair for 2D Euler Equation}

We consider the two-dimensional incompressible Euler equation
written in vorticity form,
\\
\begin{eqnarray}
& &{\pa \Om \over \pa t}= - u \ {\pa \Om \over \pa x} -v
\ {\pa \Om \over \pa y}\ ,
\nonumber \\
\label{Euler} \\
& &{\pa u \over \pa x} +{\pa v \over \pa y} = 0\ ; \nonumber
\end{eqnarray}
where $\Om$ is vorticity, $u$ and $v$ are 
respectively velocity components along $x$ and $y$ directions.
If we define the stream function $\Psi$ by,
\\
\[
u=- {\pa \Psi \over \pa y} \ ,\ \ \ v={\pa \Psi \over \pa x} \ ,
\]
\\
\nid
then we have the relation between vorticity $\Om$ and stream 
function $\Psi$,
\\
\[
\Om ={\pa v \over \pa x} -{\pa u \over \pa y} =\Dl \Psi \ .
\]
\\
\nid
The 2D Euler equation (\ref{Euler}) can be written in the equivalent 
form,
\begin{equation}
{\pa \Om \over \pa t} + \{ \Psi, \Om \} = 0 \ ,
\label{euler}
\end{equation}
where the bracket $\{\ \}$ is defined as
\[
\{ f, g\} = (\pa_x f) (\pa_y g) - (\pa_y f) (\pa_x g) \ .
\]
\begin{theorem}[\cite{Li00c}]
The Lax pair of the 2D Euler equation (\ref{euler}) is given as
\begin{equation}
\left \{ \begin{array}{l} 
L \varphi = \la \varphi \ ,
\\
\pa_t \varphi + A \varphi = 0 \ ,
\end{array} \right.
\label{laxpair}
\end{equation}
where
\[
L \varphi = \{ \Om, \varphi \}\ , \ \ \ A \varphi = \{ \Psi, \varphi \}\ ,
\]
and $\la$ is a complex constant, and $\varphi$ is a complex-valued function.
\end{theorem}
The compatibility condition of the Lax pair (\ref{laxpair}) gives the 
2D Euler equation (\ref{euler}), i.e. 
\[
\pa_t L = [L, A]\ ,
\]
where $[L, A] = LA -AL$, gives the Lax representation of the 2D Euler 
equation (\ref{euler}). In investigating 2D Euler equation through 
the Lax pair, two areas seem promising. The first area is Darboux 
transformation. Darboux transformations have been utilized for 
generating explicit representations for homoclinic structures for 
soliton equations, see for example \cite{Li00a}. Up to now, the Darboux 
transformation for the above Lax pair has not been found. On the other 
hand, in the present paper, 
we do not use the approach of Darboux transformation for investigating 
homoclinic structures for 2D Euler equation, rather use the classical 
approach of Galerkin truncations. Investigating the Darboux transformation 
is our future research.
The other area is inverse scattering 
transform. Inverse scattering transforms have been exploited for solving 
the Cauchy problems of soliton equations, see for example \cite{FA83} 
\cite{BC89}.
In building the inverse scattering, it is crucial to have constant 
coefficient differential operators in (especially the spatial part of)
the Lax pair. For that reason, one may start with the following Lax 
pair \cite{Zak89},
\begin{equation}
\left \{ \begin{array}{l} 
D_1 \varphi + \{ \Om, \varphi \} = \la \varphi \ ,
\\
\pa_t \varphi + D_2 \varphi + \{ S, \varphi \}= 0 \ ,
\end{array} \right.
\label{laxp1}
\end{equation}
where
\[
D_1 = \al {\pa \over \pa x} + \be {\pa \over \pa y} \ , \ \ \ 
D_2 = \ga {\pa \over \pa x} + \dl {\pa \over \pa y} \ ,
\] 
$\al , \be , \ga , \ \mbox{and} \ \dl $ are real constants, $\la$ is a complex 
constant, $S$ is a real-valued function, and $\varphi$ is a complex-valued 
function. The compatibility condition of this Lax pair gives the 
following equation instead of the 2D Euler equation,
\begin{equation}
{\pa \Om \over \pa t} + \{ S, \Om \}+ D_2 \Om - D_1 S = 0.
\label{veq1}
\end{equation}
A promising approach is that one can first build the inverse scattering 
for (\ref{laxp1}) and (\ref{veq1}), and then take the limits 
\[
\al , \be , \ga , \dl \ra 0\ , \ \ \ S \ra \Psi = \Dl^{-1} \Om \ ,
\]
to get results for 2D Euler equation. A special case of the Equation 
(\ref{veq1}) is the following system of equations \cite{Zak89},
\[
\left \{ \begin{array}{l} 
{\pa \Om \over \pa t} + \{ S, \Om \} = 0 \ ,
\\ \\
D_1 S = D_2 \Om\ .
\end{array} \right.
\]
In this paper, we are not going to follow the approach of inverse scattering 
transform either. Investigating the inverse scattering 
transform is our future research. Here the Lax pair 
structures make it plausible for the existence of homoclinic 
structures for the 2D Euler equation.

\subsection{The Lax Pair for 3D Euler Equation}

We consider the three-dimensional incompressible Euler equation
written in vorticity form,
\begin{equation}
\pa_t \Om + (u \cdot \na) \Om - (\Om \cdot \na) u = 0 \ ,
\label{3deuler}
\end{equation}
where $u = (u_1, u_2, u_3)$ is the velocity, $\Om = (\Om_1, \Om_2, \Om_3)$
is the vorticity, $\Om = \na \times u$, and $\na \cdot u = 0$. $u$ can be 
represented by $\Om$ for example through Biot-Savart law. More recently, 
to the author's surprise, the Lax pair for 3D Euler equation was established 
by Steve Childress.
\begin{theorem}[\cite{Chi00}]
The Lax pair of the 3D Euler equation (\ref{3deuler}) is given as
\begin{equation}
\left \{ \begin{array}{l} 
L \varphi = \la \varphi \ ,
\\
\pa_t \varphi + A \varphi = 0 \ ,
\end{array} \right.
\label{3dlaxpair}
\end{equation}
where
\[
L \varphi = \Om \cdot \na \varphi - \varphi \cdot \na \Om \ , 
\ \ \ A \varphi = u \cdot \na \varphi - \varphi \cdot \na u \ , 
\]
$\la$ is a complex constant, and $\varphi = (\varphi_1, \varphi_2, 
\varphi_3)$ is a complex 3-vector valued function.
\end{theorem}
The compatibility condition of the Lax pair (\ref{3dlaxpair}) gives the 
3D Euler equation (\ref{3deuler}), i.e. 
\[
\pa_t L = [L, A]\ ,
\]
where $[L, A] = LA -AL$, gives the Lax representation of the 3D Euler 
equation (\ref{3deuler}). 

Unfortunately, when one does 2D reduction to the Lax pair (\ref{3dlaxpair}),
one gets $L = 0$. Therefore, the Lax pair (\ref{3dlaxpair}) does not 
imply any Lax pair for 2D Euler equation.

In investigating 3D Euler equation through 
the Lax pair, especially the question on the possibility of finite time 
blow up solutions, the two areas: (1). Darboux transformations, (2). 
inverse scattering transforms, also seem promising. Up to now, the Darboux 
transformation for the Lax pair (\ref{3dlaxpair}) has not been found.
Investigating the Darboux transformation and the inverse scattering transform
for the Lax pair (\ref{3dlaxpair}) is our future research.

In building the inverse scattering, it is crucial to have constant 
coeficient differential operators in (especially the spatial part of)
the Lax pair. For 3D Euler equation, we may start with the following Lax 
pair
\begin{proposition}
We consider the following Lax pair,
\begin{equation}
\left \{ \begin{array}{l} 
L \varphi = \la \varphi \ ,
\\
\pa_t \varphi + A \varphi = 0 \ ,
\end{array} \right.
\label{3vlaxpair}
\end{equation}
where
\[
L \varphi = \Om \cdot \na \varphi - \varphi \cdot \na \Om +D_1 \varphi\ , 
\ \ \ A \varphi = q \cdot \na \varphi - \varphi \cdot \na q +D_2 \varphi\ , 
\]
$\la$ is a complex constant, $\varphi = (\varphi_1, \varphi_2, 
\varphi_3)$ is a complex 3-vector valued function, $q=(q_1, q_2, q_3)$
is a real 3-vector valued function, $D_j = \al^{(j)} \cdot \na$, $(j=1,2)$,
$\al^{(j)}= (\al^{(j)}_1, \al^{(j)}_2, \al^{(j)}_3)$ are real constant
3-vectors. The compatibility condition of this Lax pair gives the 
following equation instead of the 3D Euler equation,
\begin{equation}
\pa_t \Om + (q \cdot \na) \Om - (\Om \cdot \na) q +D_2 \Om - D_1 q = 0 \ .
\label{3veuler}
\end{equation}
A specialization of (\ref{3veuler}) is the following system of equations,
\[
\left \{ \begin{array}{l} 
\pa_t \Om + (q \cdot \na) \Om - (\Om \cdot \na) q = 0 \ ,
\\ \\
D_1 q = D_2 \Om\ .
\end{array} \right.
\]
\end{proposition}

Proof: The proof of this proposition is a trivial direct verification. $\Box$

A promising approach is that one can first build the inverse scattering 
for (\ref{3vlaxpair}), and then take the limits 
\[
\al^{(j)} \ra 0\ , \ (j=1,2)\ ,  \ \ \ q \ra u \ ,
\]
to get results for 3D Euler equation. Here the Lax pair 
structures make it plausible for the existence of homoclinic 
structures for the 3D Euler equation based upon the study on 
soliton equations \cite{Li00a}.
\begin{remark}
If one defines integrability of an equation 
by the existence of a Lax pair, then both 2D and 3D Euler equations 
are integrable. More importantly, both 2D and 3D Navier-Stokes equations 
at high Reynolds numbers are near integrable systems.
\end{remark}

\eqnsection{Introduction of the Dashed-Line Model}

We consider the two-dimensional incompressible Euler equation
written in vorticity form(\ref{Euler})
under periodic boundary conditions in both $x$ and $y$ directions
with period $2\pi$. We also require that both 
$u$ and $v$ have means zero,
\\
\[
\int_0^{2\pi}\int_0^{2\pi} u\ dxdy =\int_0^{2\pi}\int_0^{2\pi} v\ dxdy=0.
\]
\\
\nid
We expand $\Om$ into Fourier series,
\\
\[
\Om =\sum_{k\in Z^2/\{0\}} \om_k \ e^{ik\cdot X}\ ,
\]
\\
\nid
where $\om_{-k}=\overline{\om_k}\ $, $k=(k_1,k_2)^T$, 
$X=(x,y)^T$. In this paper, we confuse $0$ with $(0,0)^T$, the context 
will always make it clear. By the relation between vorticity $\Om$ and stream 
function $\Psi$,
the system (\ref{Euler}) can be rewritten as the following kinetic system,
\\
\begin{equation}
\dot{\om}_k = \sum_{k=p+q} A(p,q) \ \om_p \om_q \ ,
\label{Keuler}
\end{equation}
\\
\nid
where $A(p,q)$ is given by,
\\
\begin{eqnarray}
A(p,q)&=& {1\over 2}[|q|^{-2}-|p|^{-2}](p_1 q_2 -p_2 q_1) \nonumber \\
\label{Af} \\      
      &=& {1\over 2}[|q|^{-2}-|p|^{-2}]\left | \begin{array}{lr} 
p_1 & q_1 \\ p_2 & q_2 \\ \end{array} \right | \ , \nonumber
\end{eqnarray}
\\
\nid
where $|q|^2 =q_1^2 +q_2^2$ for $q=(q_1,q_2)^T$, similarly for $p$.

For any two functionals $F_1$ and $F_2$ of $\{ \om_k \}$, we 
define their Lie-Poisson bracket as
\\
\begin{equation}
\{ F_1,F_2 \} = \sum_{k+p+q=0} \left | \begin{array}{lr} 
q_1 & p_1 \\ q_2 & p_2 \\ \end{array} \right | \ \om_k \
{\pa F_1 \over \pa \overline{\om_p}} \ {\pa F_2 \over \pa \overline{\om_q}}\ .
\label{Liebr}
\end{equation}
\\
\nid
Then the 2D Euler equation (\ref{Keuler}) is a Hamiltonian system \cite{Arn66},
\\
\begin{equation}
\dot{\om}_k = \{ \om_k, H\}, \label{hEft}
\end{equation}
\\
\nid
where the Hamiltonian $H$ is the kinetic energy,
\\
\begin{equation}
H= {1\over 2} \sum_{k \in Z^2/\{0\}} |k|^{-2} |\om_k |^2. 
\label{Ih}
\end{equation}
\\
\nid
Following are Casimirs (i.e. invariants that Poisson commute with 
any functional) of the Hamiltonian system (\ref{hEft}):
\\
\begin{equation}
J_n = \sum_{k_1 + \cdot \cdot \cdot +k_n =0} \om_{k_1} 
\cdot \cdot \cdot \om_{k_n}. \label{Ic}
\end{equation}

\subsection{Preliminaries on Linearized 2D Euler Equation}

In this subsection, we discuss the preliminary results on linearized 
2D Euler equation known from \cite{Li00}.

We denote $\{ \om_k \}_{k\in \Z}$ by $\om$. We consider the simple fixed point 
$\om^*$:
\\
\begin{equation}
\om^*_p = \Ga,\ \ \ \om^*_k = 0 ,\ \mbox{if} \ k \neq p \ \mbox{or}\ -p,
\label{fixpt}
\end{equation}
\\
\nid
of the 2D Euler equation (\ref{Keuler}), where 
$\Ga$ is an arbitrary complex constant. 
The {\em{linearized two-dimensional Euler equation}} at $\om^*$ is given by,
\\
\begin{equation}
\dot{\om}_k = A(p,k-p)\ \Ga \ \om_{k-p} + A(-p,k+p)\ \bar{\Ga}\ \om_{k+p}\ .
\label{LE}
\end{equation}
\begin{definition}[Classes]
For any $\hk \in \Z$, we define the class $\Sg_{\hk}$ to be the subset of 
$\Z$:
\\
\[
\Sg_{\hk} = \bigg \{ \hk + n p \in \Z \ \bigg | \ n \in Z, \ \ p \ \mbox{is 
specified in (\ref{fixpt})} \bigg \}.
\]
\label{classify}
\end{definition}
See Fig.\ref{class} for an illustration of the classes. 
According to the classification 
defined in Definition \ref{classify}, the linearized two-dimensional Euler 
equation (\ref{LE}) decouples into infinitely many {\em{invariant subsystems}}:
\\
\begin{eqnarray}
\dot{\omega}_{\hat{k} + np} &=& A(p, \hat{k} + (n-1) p) 
     \ \Gamma \ \omega_{\hat{k} + (n-1) p} \nonumber \\  \label{CLE}\\
& & + \ A(-p, \hat{k} + (n+1)p)\ 
     \bar{\Gamma} \ \omega_{\hat{k} +(n+1)p}\ . \nonumber
\end{eqnarray}
\begin{figure}[ht]
  \begin{center}
    \leavevmode
      \setlength{\unitlength}{2ex}
  \begin{picture}(36,27.8)(-18,-12)
    \thinlines
\multiput(-12,-11.5)(2,0){13}{\line(0,1){23}}
\multiput(-16,-10)(0,2){11}{\line(1,0){32}}
    \thicklines
\put(0,-14){\vector(0,1){28}}
\put(-18,0){\vector(1,0){36}}
\put(0,15){\makebox(0,0){$k_2$}}
\put(18.5,0){\makebox(0,0)[l]{$k_1$}}
%
%
%
%
%
\qbezier(-5.5,0)(-5.275,5.275)(0,5.5)
\qbezier(0,5.5)(5.275,5.275)(5.5,0)
\qbezier(5.5,0)(5.275,-5.275)(0,-5.5)
\qbezier(0,-5.5)(-5.275,-5.275)(-5.5,0)
    \thinlines
\put(4,4){\circle*{0.5}}
\put(0,0){\vector(1,1){3.7}}
\put(4.35,4.35){$p$}
\put(4,-4){\circle*{0.5}}
\put(8,0){\circle*{0.5}}
\put(-8,0){\circle*{0.5}}
\put(-8,-2){\circle*{0.5}}
\put(-12,-4){\circle*{0.5}}
\put(-12,-6){\circle*{0.5}}
\put(-4,2){\circle*{0.5}}
\put(-4,4){\circle*{0.5}}
\put(0,6){\circle*{0.5}}
\put(0,8){\circle*{0.5}}
\put(4,10){\circle*{0.5}}
\put(12,4){\circle*{0.5}}
\put(0,-8){\circle*{0.5}}
\put(-4,-12){\line(1,1){17.5}}
\put(-13.5,-7.5){\line(1,1){19.5}}
\put(-13.5,-5.5){\line(1,1){17.5}}
\put(-3.6,1.3){$\hat{k}$}
\put(-7,12.1){\makebox(0,0)[b]{$(-p_2, p_1)^T$}}
\put(-6.7,12){\vector(1,-3){2.55}}
\put(6.5,13.6){\makebox(0,0)[l]{$\Sg_{\hat{k}}$}}
\put(6.4,13.5){\vector(-2,-3){2.0}}
\put(7,-12.1){\makebox(0,0)[t]{$(p_2, -p_1)^T$}}
\put(6.7,-12.25){\vector(-1,3){2.62}}
\put(-4.4,-13.6){\makebox(0,0)[r]{$\bar{D}_{|p|}$}}
\put(-4.85,-12.55){\vector(1,3){2.45}}
  \end{picture}
  \end{center}
\caption{An illustration of the classes $\Sg_{\hk}$ and the disk 
$\bar{D}_{|p|}$.}
\label{class}
\end{figure}
\begin{theorem}
The eigenvalues of the linear operator $\LL_{\hk}$ defined by the 
right hand side 
of (\ref{CLE}), are of 
four types: real pairs ($c, -c$), purely imaginary pairs ($id, -id$), 
quadruples ($\pm c \pm id$), and zero eigenvalues.
\end{theorem}
The eigenvalues can be computed through continued fractions.
\begin{definition}[The Disk]
The disk of radius $\left| p \right|$ in $Z^2 / \left\{ 0
\right\}$, denoted by $D_{\left| p \right|}$, is defined as
\[
  D_{\left| p \right|} = \bigg \{ k \in Z^2 / \left\{ 0 \right\} \ \bigg| 
      \ \left| k \right| < \left| p \right| \bigg \} \, .
\]
The closure of $D_{\left| p \right|}$, denoted by
$\bar{D}_{\left| p \right|}$, is defined as
\[ 
 \bar{D}_{\left| p \right|} = \bigg \{ k \in Z^2/ \left\{ 0 \right\} \ \bigg| 
     \ \left| k \right| \leq \left| p \right| \bigg \} \, .
\]
\end{definition}
See Fig.\ref{class} for an illustration. 
\begin{theorem}[The Spectral Theorem] We have the following claims on 
the spectra of the linear operator $\LL_{\hk}$:
\begin{enumerate}
\item If $\Sg_{\hat{k}} \cap \bar{D}_{|p|} = \emptyset$, then the entire
$\ell_2$ spectrum of the linear operator $\LL_{\hk}$ 
is its continuous spectrum. See Figure \ref{splb}, where
$b= - \frac{1}{2}|\Gamma | |p|^{-2} 
\left|
  \begin{array}{cc}
p_1 & \hat{k}_1 \\
p_2 & \hat{k}_2
  \end{array}
\right| \ .$
\item If $\Sg_{\hat{k}} \cap \bar{D}_{|p|} \neq \emptyset$, then the entire
essential $\ell_2$ spectrum of the linear operator $\LL_{\hk}$ is its 
continuous spectrum. 
That is, the residual 
spectrum of $\LL_{\hk}$ is empty, $\sg_r (\LL_{\hk}) = \emptyset$. The point 
spectrum of $\LL_{\hk}$ is symmetric with respect to both real and 
imaginary axes. 
See Figure \ref{spla2}.
\end{enumerate}
\label{spthla}
\end{theorem}
\begin{figure}[ht]
  \begin{center}
    \leavevmode
      \setlength{\unitlength}{2ex}
  \begin{picture}(36,27.8)(-18,-12)
    \thicklines
\put(0,-14){\vector(0,1){28}}
\put(-18,0){\vector(1,0){36}}
\put(0,15){\makebox(0,0){$\Im \{ \la \}$}}
\put(18.5,0){\makebox(0,0)[l]{$\Re \{ \la \}$}}
\put(0.1,-7){\line(0,1){14}}
\put(.2,-.2){\makebox(0,0)[tl]{$0$}}
\put(-0.2,-7){\line(1,0){0.4}}
\put(-0.2,7){\line(1,0){0.4}}
\put(2.0,-6.4){\makebox(0,0)[t]{$-i2|b|$}}
\put(2.0,7.6){\makebox(0,0)[t]{$i2|b|$}}
\end{picture}
  \end{center}
\caption{The spectrum of $\LL_{\hk}$ in case 1.}
\label{splb}
\end{figure}
\begin{figure}[ht]
  \begin{center}
    \leavevmode
      \setlength{\unitlength}{2ex}
  \begin{picture}(36,27.8)(-18,-12)
    \thicklines
\put(0,-14){\vector(0,1){28}}
\put(-18,0){\vector(1,0){36}}
\put(0,15){\makebox(0,0){$\Im \{ \la \}$}}
\put(18.5,0){\makebox(0,0)[l]{$\Re \{ \la \}$}}
\put(0.1,-7){\line(0,1){14}}
\put(.2,-.2){\makebox(0,0)[tl]{$0$}}
\put(-0.2,-7){\line(1,0){0.4}}
\put(-0.2,7){\line(1,0){0.4}}
\put(2.0,-6.4){\makebox(0,0)[t]{$-i2|b|$}}
\put(2.0,7.6){\makebox(0,0)[t]{$i2|b|$}}
\put(2.4,3.5){\circle*{0.5}}
\put(-2.4,3.5){\circle*{0.5}}
\put(2.4,-3.5){\circle*{0.5}}
\put(-2.4,-3.5){\circle*{0.5}}
\put(5,4){\circle*{0.5}}
\put(-5,4){\circle*{0.5}}
\put(5,-4){\circle*{0.5}}
\put(-5,-4){\circle*{0.5}}
\put(8,6){\circle*{0.5}}
\put(-8,6){\circle*{0.5}}
\put(8,-6){\circle*{0.5}}
\put(-8,-6){\circle*{0.5}}
\end{picture}
  \end{center}
\caption{The spectrum of $\LL_{\hk}$ in case 2.}
\label{spla2}
\end{figure}

\subsection{Rationality of the Dashed-Line Model}

To simplify our study, we study only the case when $\om_k$ is real, $\forall 
k \in \Z$, i.e. we only study the cosine transform of the vorticity, 
\[
\Om = \sum_{k \in \Z} \om_k \cos (k \cdot X)\ ,
\]
and the 2D Euler equation (\ref{Euler};\ref{Keuler}) preserves the cosine 
transform. To further simplify our study, we will study a concrete 
dashed-line model based upon the line of fixed points (\ref{fixpt}) with the 
mode $p=(1,1)^T$ parametrized by $\Ga$.
When $\Ga \neq 0$, each fixed point has $4$ eigenvalues which form a 
quadruple. These four eigenvalues appear in the only unstable invariant 
linear subsystem labeled by $\hk = (-3,-2)^T$. We computed the eigenvalues 
through continued fractions, one of them is \cite{Li00}:
\begin{equation}
\tla=2 \lambda / | \Gamma | = 0.24822302478255 \ + \ i \ 0.35172076526520\ .
\label{evun}
\end{equation}
See Figure \ref{figev} for an illustration.
\begin{figure}[ht]
  \begin{center}
    \leavevmode
      \setlength{\unitlength}{2ex}
  \begin{picture}(36,27.8)(-18,-12)
    \thicklines
\put(0,-14){\vector(0,1){28}}
\put(-18,0){\vector(1,0){36}}
\put(0,15){\makebox(0,0){$\Im \{ \tla \}$}}
\put(18.5,0){\makebox(0,0)[l]{$\Re \{ \tla \}$}}
\put(2.4,3.5){\circle*{0.5}}
\put(-2.4,3.5){\circle*{0.5}}
\put(2.4,-3.5){\circle*{0.5}}
\put(-2.4,-3.5){\circle*{0.5}}  
\end{picture}
  \end{center}
\caption{The quadruple of eigenvalues of the invariant system 
labeled by $\hk = (-3,-2)^T$, when $p=(1,1)^T$.}
\label{figev}
\end{figure}
We hope that a Galerkin truncation with a small 
number of modes including those inside the disk 
$\bar{D}_{\left| p \right|}$ can capture the eigenvalues.
We propose the Galerkin truncation to the linear system 
(\ref{CLE}) with the four modes $\hk +p$,  $\hk +2p$, $\hk +3p$, and $\hk +4p$,
\begin{eqnarray*}
\dot{\omega}_{1} &=& -A_{2} \Gamma \omega_{2} \, , \\
\dot{\omega}_{2} &=& A_{1} \Gamma \omega_{1}
- A_{3} \Gamma \omega_{3}\, , \\
\dot{\omega}_{3} &=& A_{2} \Gamma \omega_{2}
- A_{4} \Gamma \omega_{4}\, , \\
\dot{\omega}_{4} &=& A_{3} \Gamma \omega_{3} \, . 
\end{eqnarray*}
From now on, the abbreviated notations,
\begin{equation}
  \omega_n = \omega_{\hat{k}+np} \, , \ \  
  A_n = A(p,\hat{k}+np) \, , \ \  
  A_{m,n} = A(\hat{k}+mp,\hat{k}+np) \, ,  
\label{abbn}
\end{equation}
will be used. The eigenvalues of this four dimensional system can be 
easily calculated. It turns out that this system has a quadruple of 
eigenvalues:
\begin{eqnarray}
\lambda &=& \pm \frac{\Gamma}{2 \sqrt{10}} \sqrt{1 \pm i \sqrt{35}} 
\nonumber\\
&\dot{=}& \pm\left( \frac{\Gamma}{2} \right) \times 0.7746 
  \times e^{\pm i \theta_1} \, , \label{evu}
\end{eqnarray}
where $\theta_1 = \arctan (0.845)$, in comparison with the
quadruple of eigenvalues (\ref{evun}), where
\begin{displaymath}
  \lambda \dot{=} \pm \left( \frac{\Gamma}{2} \right) \times 0.43 \times 
  e^{\pm i \theta_2} \, ,
\end{displaymath}
and $\theta_2 = \arctan (1.418)$. Thus, {\em{ the quadruple of eigenvalues of
the original system is recovered by the four-mode truncation}}. 
Along the above line of thinking, we ``chop'' the 
line $\Sg_{\hk}$ as follows,
\begin{eqnarray*}
\dot{\omega}_{5j+1} &=& -A_{5j+2} \ \Gamma \ \omega_{5j+2} \, , \\
\dot{\omega}_{5j+2} &=& A_{5j+1} \ \Gamma \ \omega_{5j+1}
  - A_{5j+3} \ \Gamma \ \omega_{5j+3}\, , \\
\dot{\omega}_{5j+3} &=& A_{5j+2} \ \Gamma \ \omega_{5j+2}
  - A_{5j+4} \ \Gamma \ \omega_{5j+4}\, , \\
\dot{\omega}_{5j+4} &=& A_{5j+3} \ \Gamma \ \omega_{5j+3}\, . 
\end{eqnarray*}
That is, we ``chop'' the line $\Sg_{\hk}$ at the points $n=5j,\ \forall 
j \in Z$. When $j \neq 0$, this system has two complex conjugate pairs 
of purely imaginary eigenvalues,
\[
\la = \pm i (\Ga / \sqrt{2})\sqrt{\bigg | b +\dl \ \sqrt{b^2-4c} \bigg |}\ ,
\ \ \ \ \dl = \pm \ ;
\]
where
\[
b = - A_{5j+1}A_{5j+2} - A_{5j+2}A_{5j+3} - A_{5j+3}A_{5j+4} \ ,
\]
\[
c = A_{5j+1}A_{5j+2}A_{5j+3}A_{5j+4}\ ,
\]
\[
b < 0\ , \ \ c > 0\ , \ \ 0 < b^2 -4c < b^2\ .
\]
When $j \ra \pm \infty$, 
\begin{eqnarray*}
\la &\ra& \pm i (\Ga /2)\sqrt{\frac{3 +\dl \ \sqrt{5}}{2}}\ ,
\ \ \ \ \dl = \pm \ , \\
&\dot{=}& \pm \ i \ 0.31 \ \Ga\ , \ \ \ \pm \ i \ 0.8 \ \Ga \ .
\end{eqnarray*}
For example, when $j=1$,
\[
\la \ \dot{=} \ \pm \ i \ 0.2937609 \ \Ga\ , \ \  \pm \ i \ 0.7736967 \ \Ga\ ;
\]
when $j=2$,
\[
\la \ \dot{=} \ \pm \ i \ 0.3057701 \ \Ga\ , \ \  \pm \ i \ 0.8007493 \ \Ga\ .
\]
For all $j \in Z$, the distribution of the eigenvalues is illustrated 
in Fig.\ref{disev}.
\begin{figure}
\vspace{1.5in}
\caption{The distribution of the eigenvalues of the chopped system.}
\label{disev}
\end{figure}
To recover the original linear system from the ``chopped'' system, 
we introduce the homotopy parameter $\e \in [0,1]$, and consider the 
homotopy system,
\[
\dot{\omega}_n = \epsilon_{n-1} A_{n-1} \Ga
  \omega_{n-1} - \epsilon_{n+1} A_{n+1} \Ga \omega_{n+1} \,
\]
where 
\[
\epsilon_n = \left\{
    \begin{array}{ll}
      1 \, , & \ \ \hbox{if } n \ne 5j \, , \, \forall j \in Z \ ,  \\
\epsilon \, , & \ \ \hbox{if } n = 5j \, , \, 
      \hbox{ for some }j \in Z \, .
    \end{array} \right.
\]
The homotopy deformation of the spectra for this homotopy system is 
calculated numerically by Thomas Witelski and shown in Fig. \ref{sphty} 
for a special value of $\Ga$. 
\begin{figure}
\vspace{1.5in}
\caption{The homotopy deformation of the spectra for the homotopy system.}
\label{sphty}
\end{figure}
As $\e$ is increased from $0$ to $1$, the quadruple of eigenvalues 
moves along a curve, and the sizes of the four bands 
of eigenvalues with zero real parts increase, and finally these four bands 
congregate into one band when $\e = 1$
which is the continuum spectrum. From such numerical 
calculation, one can clearly see the homotopy transition of the spectra 
from a decoupled ``chopped'' subsystem to the original linearized 2D Euler 
equation. From the above discussions, it is 
natural to propose the following {\em{dashed-line model}} ,
\begin{eqnarray}
\dot{\omega}_n &=& \epsilon_{n-1} A_{n-1} \omega_p
\omega_{n-1} - \epsilon_{n+1} A_{n+1} \omega_p \omega_{n+1} \ , 
\nonumber \\ 
\label{rdlm} \\
\dot{\omega}_p &=& - \sum_{n \in Z} \epsilon_n \epsilon_{n-1}
A_{n-1,n} \omega_{n-1} \omega_n \, , \nonumber
\end{eqnarray}
to model the hyperbolic structure of the 2D Euler equation, connected to 
the line of fixed points (\ref{fixpt}) with $p=(1,1)^T$. Figure \ref{model}
illustrates the collocation of the modes in this model, which has the 
{\em{``dashed-line''}} nature leading to the name of the model. 
\begin{figure}[ht]
  \begin{center}
    \leavevmode
      \setlength{\unitlength}{2ex}
  \begin{picture}(36,27.8)(-18,-12)
    \thinlines
\multiput(-12,-11.5)(2,0){13}{\line(0,1){23}}
\multiput(-16,-10)(0,2){11}{\line(1,0){32}}
    \thicklines
\put(0,-14){\vector(0,1){28}}
\put(-18,0){\vector(1,0){36}}
\put(0,15){\makebox(0,0){$k_2$}}
\put(18.5,0){\makebox(0,0)[l]{$k_1$}}
\qbezier(-2.75,0)(-2.6375,2.6375)(0,2.75)
\qbezier(0,2.75)(2.6375,2.6375)(2.75,0)
\qbezier(2.75,0)(2.6375,-2.6375)(0,-2.75)
\qbezier(0,-2.75)(-2.6375,-2.6375)(-2.75,0)
    \thinlines
\put(2,2){\circle*{0.5}}
\put(0,0){\vector(1,1){1.85}}
\put(2.275,2.275){$p$}
\put(-12,-10){\circle*{0.5}}
\put(-10,-8){\circle*{0.5}}
\put(-8,-6){\circle*{0.5}}
\put(-6,-4){\circle{0.5}}
\put(-4,-2){\circle*{0.5}}
\put(-2,0){\circle*{0.5}}
\put(0,2){\circle*{0.5}}
\put(2,4){\circle*{0.5}}
\put(4,6){\circle{0.5}}
\put(6,8){\circle*{0.5}}
\put(8,10){\circle*{0.5}}
\put(-14,-12){\line(1,1){24}}
\put(-5.6,-5.4){$\hat{k}$}
\put(-4.4,-13.6){\makebox(0,0)[r]{$\bar{D}_{|p|}$}}
\put(-4.85,-12.55){\vector(1,3){3.4}}
\end{picture}
\end{center}
\caption{The collocation of the modes in the dashed-line model.}
\label{model}
\end{figure}

The dashed-line model has the same properties as 2D Euler equation, of 
conserving the kinetic energy and enstrophy, and being a Hamiltonian 
system with the same Lie-Poisson bracket structure. For any two functionals 
$F_1$ and $F_2$, we define their Lie-Poisson bracket as follows,
\[
\{ F_1,F_2\} = \sum_{k+q+r=0} \left | \begin{array}{lr}r_1 &q_1 \\ 
r_2 & q_2 \\ \end{array} \right | \e_k \om_k \frac {\pa F_1} {\pa \om_{-q}} 
 \frac {\pa F_2} {\pa \om_{-r}} \ .
\]
The kinetic energy is given by,
\[
\tH = \frac{1}{2} \bigg \{ \sum_{n \in Z} \e_n |\hk+np|^{-2} \om_{\hk +np} 
\om_{-[\hk +np]} +|p|^{-2} \om_p \om_{-p} \bigg \} \ .
\]
The dashed-line model (\ref{rdlm}) is a Hamiltonian system with the kinetic 
energy as the Hamiltonian,
\[
\left \{ \begin{array}{l} 
\dot{\om}_{\hk+np} = \{ \om_{\hk+np}, \tH\}\ , \\ \\
\dot{\om}_p = \{ \om_p, \tH\}\ . \\
\end{array} \right .
\]
The enstrophy 
\[
\tilde{J}_2 = \om_p \om_{-p} +\sum_{n\in Z}\e_n \om_{\hk+np}\om_{-[\hk+np]} 
\]
is still a constant of motion. The invariance of kinetic energy and 
enstrophy is a consequence of the simple relations,
\[
\left \{ \begin{array}{l} 
\frac{1}{|k|^2}A(p,q) + \frac{1}{|p|^2}A(q,k) + \frac{1}{|q|^2}A(k,p) = 0\ , \\
A(p,q) +  A(q,k) +  A(k,p) = 0\ ;\\
\end{array} \right .
\]
whenever $k+p+q = 0$.

The dashed-line model is also a homotopy system parametrized by the homotopy 
parameter $\e \in [0,1]$. When $\e =0$, the figure-eight 
hyperbolic structure of the 
dashed-line model is computable as studied in next section. As a first 
step toward understanding the homotopy deformation of such hyperbolic 
structure, we are going to use higher order Melnikov functions to study 
the persistence of such figure-eight hyperbolic structure when $\e$ is small. 

\eqnsection{The Degenerate Hyperbolic Foliations of the Dashed-Line 
Model When $\e = 0$}

When $\epsilon =0$, the dashed-line model takes the form:
\begin{eqnarray}
  \dot{\omega}_{5j+1} &=& -A_{5j+2} \omega_p \omega_{5j+2} \, ,\nonumber \\
 \dot{\omega}_{5j+2} &=& A_{5j+1} \omega_p \omega_{5j+1} 
       - A_{5j+3} \omega_p \omega_{5j+3} \, , \nonumber \\
 \dot{\omega}_{5j+3} &=& A_{5j+2} \omega_p \omega_{5j+2} 
       - A_{5j+4} \omega_p \omega_{5j+4} \, , \,  \ \ 
       (j \in Z) \, , \label{udlm} \\
 \dot{\omega}_{5j+4} &=& A_{5j+3} \omega_p \omega_{5j+3} \, , \nonumber \\ 
  \dot{\omega}_p &=& -\sum_{j \in Z} \bigg [ A_{5j+1, 5j+2}
    \omega_{5j+1} \omega_{5j+2} + A_{5j+2,5j+3}
    \omega_{5j+2} \omega_{5j+3} \nonumber  \\
&+& A_{5j+3, 5j+4} \omega_{5j+3} \omega_{5j+4} \bigg ] \, , \nonumber
\end{eqnarray}
and the equation for the decoupled variable $\omega_{5j}$ is given by,
\begin{displaymath}
  \dot{\omega}_{5j} = A_{5j-1} \omega_p \omega_{5j-1} 
  -A_{5j+1} \omega_p \omega_{5j+1} \, .
\end{displaymath}
This system has a sequence of invariant subsystems, one of which carries 
the hyperbolic structure.

\subsection{Invariant Subsystems When $\epsilon =0$}

For each fixed $j \in Z$, we have the seven dimensional invariant subsystems:
\begin{eqnarray}
\dot{\omega}_{5j+1} &=& -A_{5j+2} \omega_p \omega_{5j+2} \, ,\nonumber \\
\dot{\omega}_{5j+2} &=& A_{5j+1} \omega_p \omega_{5j+1} -
A_{5j+3} \omega_p  \omega_{5j+3}\, , \nonumber \\
\dot{\omega}_{5j+3} &=& A_{5j+2} \omega_p \omega_{5j+2} -
A_{5j+4} \omega_p  \omega_{5j+4}\, , \label{invsb} \\
\dot{\omega}_{5j+4} &=& A_{5j+3} \omega_p  \omega_{5j+3}\, ,\nonumber \\
\dot{\omega}_p &=& - [ A_{5j+1, 5j+2} \omega_{5j+1}
\omega_{5j+2} +  A_{5j+2, 5j+3} \omega_{5j+2}
\omega_{5j+3} \nonumber \\
&+&  A_{5j+3, 5j+4} \omega_{5j+3}
\omega_{5j+4} ] \, , \nonumber
\end{eqnarray}
and the equations for the decoupled variables $\omega_{5j}$ and
$\omega_{5j+5}$ are given by,
\begin{eqnarray*}
  \dot{\omega}_{5j} &=& - A_{5j+1} \omega_p \omega_{5j+1} \, , \\
  \dot{\omega}_{5j+5} &=& A_{5j+4}  \omega_p \omega_{5j+4} \, .
\end{eqnarray*}
The kinetic energy
\[
\tH_j = \frac{1}{2} \bigg [ \sum_{1\leq l \leq 4}\frac{1}{|\hk+(5j+l)p|^2} 
\om^2_{5j+l} + \frac{1}{|p|^2} 
\om^2_p \bigg ]\ ,
\]
and the enstrophy 
\[
\tilde{J}_2^{(j)} = \sum_{1\leq l \leq 4}\om^2_{5j+l} + \om^2_p \ ,
\]
are constants of motion for this system.

\subsection{Perverted Heteroclinic Orbits When $\epsilon =0$}

From the studies in last subsection, the ($j=0$) linearized system at 
the fixed point (\ref{fixpt}) has a quadruple of eigenvalues with nonzero 
real parts, and all ($j \neq 0$) linearized system at the fixed point 
(\ref{fixpt}) has two pairs of purely imaginary eigenvalues.
Next we focus on the invariant subsystem (\ref{invsb}) when
$j=0$, and generate explicit expressions for the unstable
and stable manifolds of the line of fixed points (\ref{fixpt}).  Notice that
\begin{eqnarray}
& & A_1 = - \frac{3}{10}\ , \  A_2 = \frac{1}{2} \ , \ 
  A_3 = A_2 \ , \ A_4 = A_1 \ ,  \nonumber \\ \label{efu} \\
& & A_{1,2} = A_1 -A_2 = - \frac{4}{5} \ , \ 
  A_{2,3} =0 \ , \ A_{3,4} = -A_{1,2} \ ; \nonumber 
\end{eqnarray}
then the invariant subsystem (\ref{invsb}) for $j=0$ can be
rewritten as follows:
\begin{eqnarray}
\dot{\omega}_1 &=& -A_2 \ \omega_p \ \omega_2 \ , \nonumber\\
\dot{\omega}_2 &=& A_1 \ \omega_p \ \omega_1 -A_2 \ \omega_p \ \omega_3
\ , \nonumber \\
\dot{\omega}_3 &=& A_2 \ \omega_p \ \omega_2 -A_1 \ \omega_p \ \omega_4
\, , \label{invu} \\
\dot{\omega}_4 &=& A_2 \ \omega_p  \ \omega_3 \, , \nonumber \\
\dot{\omega}_p &=& A_{1,2} \ (\omega_3 \ \omega_4 - \omega_1
\ \omega_2) \, , \nonumber 
\end{eqnarray}
and the equations for the decoupled variables $\omega_0$ and $\omega_5$
are given by,
\begin{eqnarray*}
\dot{\omega}_0 &=& -A_1 \ \omega_p \ \omega_1 \, , \\
\dot{\omega}_5 &=& A_1 \ \omega_p \ \omega_4 \, .
\end{eqnarray*}
There are three invariants for the system (\ref{invu}):
\begin{eqnarray}
I &=& 2 A_{1,2} (\omega_1 \omega_3 + \omega_2 \omega_4 )
        + A_2 \omega^2_p \, , \label{Iinv} \\[1ex]
U &=&  A_1 (\omega^2_1+ \omega^2_4 ) 
        + A_2 ( \omega^2_2 + \omega^2_3 )  \, , \label{Uinv} \\[1ex]
J &=&   \omega^2_p+ \omega^2_1  
        + \omega^2_2 + \omega^2_3 + \omega^2_4  \, . \label{Jinv}
\end{eqnarray}
$J$ is the enstrophy, and $U$ is a linear combination of the kinetic 
energy and the enstrophy. $I$ is an extra invariant which is peculiar 
to this invariant subsystem. With $I$, the explicit formula for the 
hyperbolic structure can be computed.

The common level set of these three invariants which is connected 
to the fixed point (\ref{fixpt}) determines the stable and
unstable manifolds of the fixed point and its negative
$-\omega^*$:
\begin{equation}
\omega_p = - \Gamma \, , \, \omega_n=0 \quad (n \in Z) \, .
\label{nfixpt}
\end{equation}
Using the polar coordinates:
\begin{displaymath}
  \omega_1 = r \cos \theta \, , \, 
  \omega_4 = r \sin \theta \, ; \, 
  \omega_2 = \rho \cos \vartheta \, , \, 
  \omega_3 = \rho \sin \vartheta \, ;
\end{displaymath}
we have the following explicit expressions for the stable and
unstable manifolds of the fixed point (\ref{fixpt}) and its
negative (\ref{nfixpt}) represented through {\em{perverted 
heteroclinic orbits}}:
\begin{eqnarray}
\omega_p &=& \Gamma \ \tanh \tau \, , \nonumber \\
r &=& \sqrt{ \frac{A_2}{A_2-A_1}}\, \ \Gamma \ \sech \tau \, , \nonumber  \\[1ex]
\theta &=& - \ \frac{A_2}{2\k} \ \mbox{ln}\ \cosh \tau + \theta_0 \, , 
\label{exus}\\[1ex]
\rho &=& \sqrt{\frac{-A_1}{A_2}} \  r \, ,\nonumber  \\[1ex]
\theta + \vartheta &=& \left\{
  \begin{array}{ll}
    - \arc \sin \left[ \frac{1}{2} \sqrt{\frac{A_2}{-A_1}}\,  \right] 
        \ , & (\k>0) \ , \\[2ex]
\pi + \arc \sin  \left[ \frac{1}{2} \sqrt{\frac{A_2}{-A_1}} \right]
\, , & (\k<0) \, ,
  \end{array} \right.  \nonumber 
\end{eqnarray}
where $A_1$ and $A_2$ are given in (\ref{efu}), $\tau = \k \Gamma t 
+ \tau_0$, $(\tau_0, \theta_0)$ are the two parameters
parametrizing the two-dimensional stable (unstable) manifold,
and
\begin{displaymath}
  \k = \sqrt{-A_1 A_2} \cos (\theta + \vartheta) 
   = \pm \sqrt{-A_1 A_2} \sqrt{1+ \frac{A_2}{4A_1}} \ .
\end{displaymath}
The two auxilliary variables $\om_0$ and $\om_5$ have the
expressions:
\begin{eqnarray*}
  \omega_0 &=& \frac{\alpha \beta}{1+ \beta^2} \sech \tau \left\{ 
            \sin [ \beta \ \mbox{ln}\ \cosh \tau +  \theta_0]
            - \frac{1}{\beta} \cos
            [ \beta \ \mbox{ln}\ \cosh \tau + \theta_0 ] \right\} \,
   , \\
  \omega_5 &=& \frac{\alpha \beta}{1+ \beta^2} \sech \tau
  \left\{ \cos [ \beta \ \mbox{ln}\ \cosh \tau + \theta_0 ]
    + \frac{1}{\beta} \sin [ \beta \ \mbox{ln}\ \cosh \tau + \theta_0]
    \right\} \, ,
\end{eqnarray*}
where
\begin{displaymath}
  \alpha = -A_1 \Gamma \k^{-1} \sqrt{\frac{A_2}{A_2-A_1}} \ , \ \ 
  \beta = - \frac{A_2}{2\k} \ .
\end{displaymath}
\begin{remark}
In this remark, we will discuss the complex version of the system 
(\ref{invu}):
\begin{eqnarray}
\dot{\omega}_1 &=& -A_2 \ \overline{\omega_p} \ \omega_2 \ , \nonumber\\
\dot{\omega}_2 &=& A_1 \ \omega_p \ \omega_1 -A_2 \ \overline{\omega_p} \ 
\omega_3 \ , \nonumber \\
\dot{\omega}_3 &=& A_2 \ \omega_p \ \omega_2 -A_1 \ \overline{\omega_p} \ 
\omega_4\, , \label{cpliv} \\
\dot{\omega}_4 &=& A_2 \ \omega_p  \ \omega_3 \, , \nonumber \\
\dot{\omega}_p &=& A_{1,2} \ (\overline{\omega_3} \ \omega_4 - 
\overline{\omega_1}\ \omega_2) \, . \nonumber 
\end{eqnarray}
This is a ten dimensional Hamiltonian system with no reality restriction. 
For this system (\ref{cpliv}), we have one complex invariant
\[
I_c = 2 A_{1,2} (\overline{\omega_1} \omega_3 + \overline{\omega_2} \omega_4 )
        + A_2 \omega^2_p \, , 
\]
and two real invariants
\begin{eqnarray*}
U_c &=&  A_1 (|\omega_1|^2+ |\omega_4|^2) 
        + A_2 (|\omega_2|^2 + |\omega_3|^2)  \, , \\[1ex]
J_c &=&   |\omega_p|^2+ |\omega_1|^2  
        + |\omega_2|^2 + |\omega_3|^2 + |\omega_4|^2  \, . 
\end{eqnarray*}
Thus, we have total four invariants which are not enough for the 
integrability of the Hamiltonian system (\ref{cpliv}). When we consider 
the reality restriction, we have three invariants. Their common level sets 
in the five dimensional real phase space, are two-dimensional. The stable 
and unstable manifolds of the fixed point are also two-dimensional, and 
these three invariants are enough to determine the stable and unstable 
manifolds.
\end{remark}

\subsection{Discussion on the Perverted Heteroclinic Orbits}

The oscillatory nature of the heteroclinic orbit reflects the fact that 
the quadruple of eigenvalues with nonzero real parts, also has nonzero 
imaginary parts. The interesting nature of the heteroclinic orbits 
is that they all have a ``perversion'' which refers to the neighborhood 
of the turning point. The portions of the heteroclinic orbit in the 
neighborhoods of $\om^*$ and $-\om^*$ can be viewed as oriented helices.
Then the entire heteroclinic orbit is a connection between a right-handed 
helix and a left-handed helix, and the connection part has to be a 
perversion. See Figure \ref{connection} for an illustration.
From the study on the linearized 2D Euler equation \cite{Li00}, we 
realize that the eigenvalues of the fixed point $\om^*$ only depend upon 
the modulus $|\om_p|$ (the same fact is true for the dashed-line model 
too). Therefore, $\om^*$ and $-\om^*$ have the same eigenvalues. In 
fact, the eigenvalues appear in quadruples which make the connection 
of a right-handed helix with a left-handed helix feasible. 

The term ``perversion'' was given by the 19th century topologist Listing 
to describe the spontaneous switching of a helical structure of one 
handedness to its mirror image. Tendril perversion in climbing plants 
had been described at length by Charles Darwin in his book \cite{Dar75}. 
Darwin interpreted the tendril perversion using the language of elasticity.
Complete mathematical study on elastic filaments modelling tendril perversion 
has been developed by A. Goriely and M. Tabor \cite{GT98a} \cite{GT96} 
\cite{GT97a} \cite{GT97b} \cite{GT97c} \cite{GT98b} \cite{GT98c}. Figure 
\ref{tendril1} shows a climbing plant with one tendril perversion 
\cite{Sac75}. One tendril can have many perversions \cite{Dar75}, see 
Figure \ref{tendril2}. Perversions can occur in problems like the microscopic 
properties of biological fibers such as cotton or the formation of bacterial 
macrofibers, telephone cords (Figure \ref{tendril1}), false-twist technique 
in the textile industry 
etc. \cite{GT98a}.
\begin{figure}
\vspace{1.5in}
\caption{The connection part of a right-handed helix and a left-handed 
helix has to be a perversion.}
\label{connection}
\end{figure}
\begin{figure}
\vspace{1.5in}
\caption{(a). Tendril perversion in Bryonia dioica. Illustration from 
Sachs' Text-book of Botany (1875); (b). Perversion in a telephone cord 
(Goriely and Tabor, 1997). The effect is achieved by fully stretching out 
and untwisting the cord and slowly bringing the ends together. The above 
figures are taken from Goriely and Tabor's papers (see References).}
\label{tendril1}
\end{figure}
\begin{figure}
\vspace{1.5in}
\caption{A caught tendril of Bryonia dioica, spirally contracted in 
reversed directions. Illustration from Darwin's book: The Movements 
and Habits of Climbing Plants (1875) (see References).}
\label{tendril2}
\end{figure}
The explicit expression (\ref{exus}) of the perverted heteroclinic orbits
shows that the unstable manifold $W^u(\om^*)$ 
of the fixed 
point $\om^*$ (\ref{fixpt}) is the same as the stable 
manifold $W^s(\om^*)$ of the negative $-\om^*$ (\ref{nfixpt}) of 
$\om^*$, and the 
stable manifold of the fixed point $\om^*$ (\ref{fixpt}) 
is the same as the unstable manifold of the negative $-\om^*$ 
(\ref{nfixpt}) of $\om^*$. Both $W^u(\om^*)$ and $W^s(\om^*)$ have 
the shapes of ``painted eggs'' (Figure \ref{het}). $W^u(\om^*)$ and 
$W^s(\om^*)$ together form the ``lip'' (Figure \ref{het}) which is a 
higher dimensional generalization of the heteroclinic connection on plane.
\begin{figure}
\vspace{1.5in}
\caption{The unstable and stable manifolds of the fixed point $\om^*$ 
(\ref{fixpt}): (a) and (b) show two ``painted eggs'' and (c) shows a 
``lip''.}
\label{het}
\end{figure}
As the first step toward understanding the degeneracy v.s. nondegeneracy 
nature of the hyperbolic structure of 2D Euler 
equation, we are interested in the $\e$-homotopy deformation of such 
hyperbolic structure for the dashed-line model. In the next section, 
we will use higher order Melnikov functions to study the ``breaking'' 
or ``persistence'' of such hyperbolic structure when $\e$ is small.

\eqnsection{The Higher Order Melnikov Functions}

In this section, we are going to study the persistence of the
hyperbolic structures given by (\ref{exus}), i.e.,~the
persistence of the ``painted eggs'' and ``lips'' as shown in
Figure~\ref{het} of the invariant subsystem (\ref{invu}) in the
dashed-line model (\ref{rdlm}) when $\epsilon \ne 0$.  We will
use Melnikov functions to detect such persistence, and the
Melnikov functions will be built upon the following two
invariants of the dashed-line model (\ref{rdlm}) when $\epsilon
=0$, which are actually invariants of the invariant subsystem
(\ref{invu}):
\begin{eqnarray}
  U &=& A_1 (\omega^2_1 + \omega^2_4) 
  + A_2 (\omega^2_2 + \omega^2_3 ) \, , \label{Umel} \\
  V &=& A_2 (\omega^2_1 + \omega^2_2 + \omega^2_3 + \omega^2_4)\nonumber \\
 & &  -2 A_{1,2} ( \omega_1 \omega_3 + \omega_2 \omega_4) \, , \label{Vmel}
\end{eqnarray}
where $V$ is a linear combination of the invariants $I$ and $J$
given in (\ref{Iinv}) and (\ref{Jinv}), $V=A_2 J-I$.

From the spectral distribution Figure \ref{disev}, when
$\epsilon =0$, the fixed point (\ref{fixpt}) or (\ref{nfixpt})
has codimension $2$ center-unstable and center-stable manifolds
$W^{cu}( \pm \Gamma)$ and $W^{cs}(\pm \Gamma)$ under the flow
(\ref{udlm}), and the ``painted eggs'' are two-dimensional
submanifolds of them.  For sufficiently small $\epsilon$, 
$W^{cu}(\pm \Gamma)$ and $W^{cs}(\pm \Gamma)$ perturb into 
$W^{cu}_\epsilon (\pm \Gamma)$ and $W^{cs}_\epsilon (\pm \Gamma)$.  We 
are going to use the new notation
\begin{equation}
  \omega = (\omega_p , \ \omega_n : n \in Z) \, .
\label{nn}
\end{equation}
In this notation, we rewrite the dashed-line model (\ref{rdlm})
as follows:
\begin{equation}
  \dot{\omega} = f(\omega) + \epsilon g(\omega) \, .
 \label{ndlm}
\end{equation}
Denote by $\omega^{(0)}$ the orbit given by (\ref{exus}) which
solves the $\epsilon =0$ form of (\ref{ndlm}):
\begin{equation}
  \dot{\omega}^{(0)} = f(\omega^{(0)}) \, .
\label{heteq}
\end{equation}

\subsection{The Melnikov Functions}

The (first order) Melnikov functions are given by:
\begin{eqnarray}
  M_U &=& \int^{\infty}_{- \infty} 
  \langle \na U , g \rangle |_{\omega^{(0)}} \, dt \label{fmlU} \\[1ex]
  &=& \int^{\infty}_{- \infty} 
  [2A_0A_1 \omega_p \omega_0 \omega_1 -
  2A_0A_1 \omega_p \omega_4 \omega_5 ]  |_{\omega^{(0)}} \, dt \nonumber 
\\[1ex]
  &=& -2A_0  \int^{\infty}_{- \infty} 
  [\omega^{(0)}_0 \dot{\omega}^{(0)}_0 
  + \omega^{(0)}_5 \dot{\omega}^{(0)}_5 ] \, dt =0 \, , \nonumber \\[1ex]
  M_V &=&  \int^{\infty}_{- \infty} 
  \langle \na V , g \rangle |_{\omega^{(0)}} \, dt \nonumber \\[1ex]
 &=&  \int^{\infty}_{- \infty} 
 [2A_0A_2 \omega_p \omega_0 \omega_1 -
 2A_0A_2 \omega_p \omega_4 \omega_5  \nonumber \\[1ex]
& & -2A_0A_{1,2} \omega_p \omega_0 \omega_3
 + 2A_0A_{1,2} \omega_p \omega_2 \omega_5 ] |_{\omega^{(0)}} \, dt \nonumber 
\\[1ex]
 &=& -2A_0  \int^{\infty}_{- \infty} \bigg [ \frac{A_2}{A_1}
   (\omega^{(0)}_0 \dot{\omega}^{(0)}_0 +
   \omega^{(0)}_5 \dot{\omega}^{(0)}_5) \nonumber \\[1ex]
& & + A_{1,2} (\omega^{(0)}_p  \omega^{(0)}_0  \omega^{(0)}_3
   - \omega^{(0)}_p  \omega^{(0)}_2  \omega^{(0)}_5) \bigg ] dt = 0 \ ,
 \nonumber 
\end{eqnarray}
since the second round bracket is an odd function in $\tau$.  The 
above representations will be reproduced in the process of
deriving higher order Melnikov functions in the next subsection.
The above calculations show that the Melnikov functions are
identically zero as functions of $\theta_0$ and $\Gamma$.  Thus
the separation between $W^{cu}_\epsilon (\pm \Gamma)$ and
$W^{cs}_\epsilon (\pm \Gamma)$ is at least of order
$O(\epsilon^2)$.  To further detect the separation between
$W^{cu}_\epsilon (\pm \Gamma)$ and  $W^{cs}_\epsilon (\pm
\Gamma)$, we need to compute the second order Melnikov functions
which are the leading order terms of the signed separation
distances. 

\subsection{The Derivation of Higher Order Melnikov Functions}

The derivation in this subsection is given without rigorous 
justifications, and rigorous justifications are future works.

Let $\omega^{(0)}(t)$ denote a heteroclinic orbit given by
(\ref{exus}), and let $\Sigma$ be a codimension $1$ hypersurface
which is transversal to the orbit  $\omega^{(0)}(t)$ at
$\omega^{(0)}(0)$.  See Figure \ref{setmel} on the setup of
this derivation.  Without loss of generality, we assume that in
backward time  $\omega^{(0)}(t)$ approaches $\omega^*$
(\ref{fixpt}) and in forward time  $\omega^{(0)}(t)$
approaches $-\omega^*$ (\ref{nfixpt}).  The submanifolds
$W^{cs}(- \Gamma) \cap \Sigma$ and  $W^{cs}_\epsilon (-
  \Gamma) \cap \Sigma$ have codimension $2$ in $\Sigma$, and the 
codimensions are coordinated by $\na U$ and
$\na V$.  Let $\omega^{(u)}(t, \epsilon)$ denote an
orbit of the $(\epsilon \ne 0)$ system (\ref{rdlm}) which
approaches $\omega^*$ in backward time.
Let $\omega^{(u)}(0, \epsilon)$ be the
intersection point of  $\omega^{(u)}(t, \epsilon)$ with
$\Sigma$.  On the submanifold $W^{cs}_\epsilon (- \Gamma) \cap
\Sigma$, let  $\omega^{(s)}(0, \epsilon)$ be the point which
has the same coordinates with $\omega^{(u)}(0, \epsilon)$
except the  $\na U$ and
$\na V$ directional coordinates in $\Sigma$.  Denote by 
$\omega^{(s)}(t, \epsilon)$ the orbit of the $(\epsilon \ne
0)$ system (\ref{rdlm}) with initial point $\omega^{(s)}(0, \epsilon)$.
\begin{figure}
\vspace{1.5in}
\caption{The geometrical setup for deriving the second order Melnikov 
functions.}
\label{setmel}
\end{figure}

Define the signed distances:
\begin{eqnarray}
  d_U &=& \langle \na U(0)  \, , \, 
         \omega^{(u)}(0, \epsilon)
         - \omega^{(s)}(0, \epsilon) \rangle \nonumber \\
      &=& \langle  \na U(0) ,
         \omega^{(u)}(0, \epsilon) - \omega^{(0)}(0)
         \rangle \nonumber \\
& & - \langle \na U(0) \, , \, 
         \omega^{(s)}(0, \epsilon)-  \omega^{(0)}(0)
         \rangle \, ,  \label{disU} \\
 d_V &=& \langle \na V(0)  \, , \, 
         \omega^{(u)}(0, \epsilon)
         - \omega^{(s)}(0, \epsilon) \rangle \nonumber \\
      &=& \langle  \na V(0) ,
         \omega^{(u)}(0, \epsilon) - \omega^{(0)}(0)
         \rangle  \nonumber \\
& & - \langle \na V(0) \, , \, 
         \omega^{(s)}(0, \epsilon)-  \omega^{(0)}(0)
         \rangle \, . \label{disV}
\end{eqnarray}
From now on, we will do the derivation for $d_U$, and the
derivation for $d_V$ is the same.  We define the
following two functions:
\begin{eqnarray*}
  \Delta^+_U(t, \epsilon) 
       &=& \langle \na U(t) ,
       \omega^{(s)}(t, \epsilon) - \omega^{(0)}(t) \rangle \, , \ \  
       t \in [0, \infty) \, , \\
   \Delta^-_U(t, \epsilon) 
       &=& \langle \na U(t, \epsilon) ,
              \omega^{(u)}(t, \epsilon)  
              -  \omega^{(0)}(t) \rangle \, , \ \  
       t \in (- \infty, 0] \, , 
\end{eqnarray*}
where $ \na U(t) = \na U(\omega^{(0)}(t))$.  Then
\begin{equation}
  d_U = \Delta^-_U(0, \epsilon) -  \Delta^+_U (0, \epsilon)
  \, .
\label{distc}
\end{equation}
Next we will do the derivation for $\Delta^+_U (t, \epsilon)$ and 
the derivation for $\Delta^-_U (t, \epsilon)$ is the
same.  For any $T>0$, $\Delta^+_U (t, \epsilon)$ as a function of
$\epsilon$ has the Taylor expansion:
\begin{equation}
  \Delta^+_U (t, \epsilon) = \sum_{n=1}^\infty \e^n \Delta^{(+,n)}_U(t)\ , 
\ \   t \in [0,T] \, ,
\label{tyU}
\end{equation}
where
\begin{equation}
    \Delta^{(+,n)}_U (t) = \frac{1}{n!}\frac{\partial^n}{\partial \epsilon^n}
    \Delta^+_U(t,0) =\frac{1}{n!} \langle \na U(t) \, , \, 
    \frac{\partial^n}{\partial \epsilon^n} \omega^{(s)}(t,0)
    \rangle \, ;
\label{dn1}
\end{equation}
and $\omega^{(s)}(t, \epsilon)$ as a function of $\epsilon$
has the Taylor expansion:
\begin{equation}
\omega^{(s)}(t, \epsilon) = \omega^{(0)}(t) + \sum_{n=1}^\infty  
\e^n \om^{(s,n)}(t)\ , \ \ t \in [0,T] \, ,
\label{tys}
\end{equation}
where
\begin{equation}
   \omega^{(s,n)}(t) =\frac{1}{n!} \frac{\partial^n}{\partial \epsilon^n} 
   \omega^{(s)}(t,0)\ .
\label{dn2}
\end{equation} 
The function $ \omega^{(s,n)}$ satisfies the equation:
\begin{equation}
\dot{\omega}^{(s,n)} = \na f(\omega^{(0)})
  \circ \omega^{(s,n)} + h^{(s,n)}\ , \ \ t \in [0,T] \, ,
\label{stbeq}
\end{equation}
where
\[
h^{(s,n)} = \frac{1}{n!} \bigg ( \frac{\pa^n}{\pa \e^n} \bigg [ 
f(\om^{(s)}(t,\e)) +\e g(\om^{(s)}(t,\e)) \bigg ] \bigg )_{\e =0} - 
\na f(\omega^{(0)})\circ \omega^{(s,n)}\ ,
\]
for example
\begin{eqnarray*}
h^{(s,1)} &=& g(\om^{(0)})\ , \\
h^{(s,2)} &=& \frac{1}{2} \bigg \{ \na^2 f(\omega^{(0)}) \circ \om^{(s,1)} 
\bigg \} \circ  \om^{(s,1)} \\
& & +\na g(\omega^{(0)}) \circ \om^{(s,1)} \ , \\
h^{(s,3)} &=& \bigg \{ \na^2 f(\omega^{(0)}) \circ \om^{(s,1)} \bigg \} 
 \circ  \om^{(s,2)}  \\
& & +\frac{1}{6} \bigg \{ \bigg \{  \na^3 f(\omega^{(0)}) \circ \om^{(s,1)}
\bigg \} \circ \om^{(s,1)} \bigg \} \circ \om^{(s,1)} \\
& & + \na g(\omega^{(0)}) \circ \om^{(s,2)} \\
& &+\frac{1}{2} \bigg \{   \na^2 g(\omega^{(0)}) \circ \om^{(s,1)} \bigg \} 
\circ \om^{(s,1)}\ .
\end{eqnarray*}

Next we are going to derive an expression for the time derivative
of the functions $\Delta^{(+,n)}_U$.
\begin{eqnarray}
& & \dot{\Delta}^{(+,n)}_U = \langle \na^2U \circ
  \dot{\omega}^{(0)} \, , \, \omega^{(s,n)} \rangle +
  \langle \na U \, , \, \dot{\omega}^{(s,n)} \rangle \nonumber \\
  &=& \langle  \na^2U \circ f \, , \, 
      \omega^{(s,n)} \rangle + \langle \na U \, , \, 
       \na f \circ   \omega^{(s,n)} \rangle
       + \langle  \na U \ , \ h^{(s,n)} \rangle \, . \label{derU}
\end{eqnarray}
Notice that
\begin{equation}
  \label{inveq}
  \langle  \na U \, , \, f \rangle =0 \, ,
\end{equation}
from the invariance of $U$ under the $\epsilon =0$ flow
(\ref{udlm}); then the variation of the equation (\ref{inveq})
leads to:
\begin{displaymath}
  \langle \na^2 U \circ \delta \omega \, , \, f \rangle
  + \langle \na U \, , \, \na f \circ 
  \delta \omega \rangle =0 \, ,
\end{displaymath}
and since $\na^2 U$ is a symmetric operator, we have
\begin{equation}
  \label{vari}
  \langle \na^2 U \circ f \, , \, \delta \omega \rangle
  + \langle \na U \, , \, \na f \circ
  \delta \omega \rangle =0 \, .
\end{equation}
If we substitute $\delta \omega$ in (\ref{vari}) by
$\omega^{(s,n)}$, we see that equation (\ref{derU}) is reduced to
\begin{equation}
  \dot{\Delta}^{(+,n)}_U = \langle \na U \, , \, h^{(s,n)} 
  \rangle \, , \, \ \ t \in [0,T] \, .
\label{derf}
\end{equation}
Without rigorous justifications and taking the limit $T \to +
\infty$ in (\ref{derf}), we have
\begin{equation}
  \Delta^{(+,n)}_U(0) = - \int^{\infty}_0 \langle 
 \na U(\omega^{(0)}) \, , \, h^{(s,n)} \rangle \, dt \, .
\label{pln}
\end{equation}
Similarly,
\begin{equation}
  \Delta^{(-,n)}_U(0) = \int_{-\infty}^0 \langle 
 \na U(\omega^{(0)}) \, , \, h^{(u,n)} \rangle \, dt \, ,
\label{mln}
\end{equation}
where $h^{(u,n)}$ has the same expression as $h^{(s,n)}$ with $s$ 
replaced by $u$ everywhere,
\begin{eqnarray}
  \Delta^-_U (t, \epsilon) &=& \sum_{n=1}^\infty \e^n \Delta^{(-,n)}_U(t)\ , 
\ \   t \in (-\infty , 0] \, , \label{dn3}  \\
\Delta^{(-,n)}_U (t) &=& \frac{1}{n!}\frac{\partial^n}{\partial \epsilon^n}
    \Delta^-_U(t,0) \nonumber \\
&=& \frac{1}{n!} \langle \na U(t) \, , \, 
    \frac{\partial^n}{\partial \epsilon^n} \omega^{(u)}(t,0)
    \rangle \, , \label{dn4} \\
\omega^{(u)}(t, \epsilon) &=& \omega^{(0)}(t) + \sum_{n=1}^\infty  
\e^n \om^{(u,n)}(t)\ , \ \ t \in (-\infty , 0] \, ,\label{dn5} \\
   \omega^{(u,n)}(t) &=&\frac{1}{n!} \frac{\partial^n}{\partial \epsilon^n} 
   \omega^{(u)}(t,0)\ ,\label{dn6} \\
\dot{\omega}^{(u,n)} &=& \na f(\omega^{(0)})
  \circ \omega^{(u,n)} + h^{(u,n)}\ , \ \ t \in  (-\infty , 0] \, .
\label{utbeq} 
\end{eqnarray}
From the expression (\ref{distc}) of the signed distance $d_U$,
we have
\begin{eqnarray*}
  d_U &=&\sum_{n=1}^\infty \epsilon^n \bigg [ \Delta^{(-,n)}_U(0)
          -  \Delta^{(+,n)}_U(0) \bigg ] \\
       &=&\sum_{n=1}^\infty \epsilon^n M^{(n)}_U\ ,
\end{eqnarray*}
where $M^{(n)}_U$ is called the n-th order Melnikov function and $M^{(1)}_U$ 
is the usual Melnikov function.  

Notice that $\om^{(u)}(0,\e)$ and $\om^{(s)}(0,\e)$ are chosen to have 
the same coordinates except the $\na U(0)$ and $\na V(0)$ directions. From 
the expressions (\ref{tyU}) (\ref{dn1}) (\ref{tys}) (\ref{dn2}) (\ref{dn3})
(\ref{dn4}) (\ref{dn5}) and (\ref{dn6}), if $M^{(j)}_Y \equiv 0$, 
$j=1,2, \cdots, n-1$; $Y=U,V$; then
\[
\om^{(s,j)}(0)=\om^{(u,j)}(0)\ , \ \ j=1,2, \cdots, n-1\ .
\]
Notice that both $\om^{(s,j)}$ and $\om^{(u,j)}$ satisfy the same form of 
equations (\ref{stbeq}) and (\ref{utbeq}); thus $\om^{(s,j)}$ and 
$\om^{(u,j)}$ together represent one solution on the entire interval 
$t \in (-\infty,\infty)$, and we can drop `$s$' and `$u$', and use 
the simple notation $\om^{(j)}$ for $j=1,2, \cdots, n-1$.
\begin{eqnarray}
\dot{\omega}^{(j)} &=& \na f(\omega^{(0)})
  \circ \omega^{(j)} + h^{(j)}\ , \ \ t \in (-\infty, \infty) \, , \label{dif}
\\
\nonumber \\
& & \ \ j=1,2, \cdots, n-1 \ , \nonumber 
\end{eqnarray}
where
\begin{eqnarray}
h^{(j)}&=& \frac{1}{j!} \bigg ( \frac{\pa^j}{\pa \e^j} \bigg [ 
f(\om^{(s)}(t,\e)) +\e g(\om^{(s)}(t,\e)) \bigg ] \bigg )_{\e =0} 
 \nonumber \\
& &- \na f(\omega^{(0)})\circ \omega^{(j)}\ , \label{hva}\\
\nonumber \\
& & \ \ j=1,2, \cdots, n \ , \nonumber
\end{eqnarray}
for example
\begin{eqnarray*}
h^{(1)} &=& g(\om^{(0)})\ , \\
h^{(2)} &=& \frac{1}{2} \bigg \{ \na^2 f(\omega^{(0)}) \circ \om^{(1)} 
\bigg \} \circ  \om^{(1)} \\
& & +\na g(\omega^{(0)}) \circ \om^{(1)} \ , \\
h^{(3)} &=& \bigg \{ \na^2 f(\omega^{(0)}) \circ \om^{(1)} \bigg \} 
 \circ  \om^{(2)}  \\
& & +\frac{1}{6} \bigg \{ \bigg \{  \na^3 f(\omega^{(0)}) \circ \om^{(1)}
\bigg \} \circ \om^{(1)} \bigg \} \circ \om^{(1)} \\
& & + \na g(\omega^{(0)}) \circ \om^{(2)} \\
& &+\frac{1}{2} \bigg \{ \na^2 g(\omega^{(0)}) \circ \om^{(1)} \bigg \} 
\circ \om^{(1)}\ .
\end{eqnarray*}
Notice that both $\om^{(u)}(t,\e)$ and $\om^{(0)}(t)$ approach the 
same point $-\om^*$ in backward time, we have
\begin{equation}
\om^{(u,l)}(-\infty) = 0 \ , \ \ \forall l = 1, 2, \cdots\ .
\label{difbc}
\end{equation}
From the expressions (\ref{pln}) and (\ref{mln}), we have
\[
M^{(n)}_U = \int^{\infty}_{-\infty} \langle \na U(\omega^{(0)}) \, , 
\, h^{(n)} \rangle \, dt \, ,
\]
similarly for $M^{(n)}_V$. 

\textbf{Summary.} {\em {The signed distances $d_Y$ $(Y=U,V)$ defined
in (\ref{disU}) and (\ref{disV}) have the representations:
\begin{equation}
  d_Y = \sum_{n=1}^\infty \epsilon^n M^{(n)}_Y \ , \quad (Y=U,V)\ .
\label{hmelf}
\end{equation}
If $M^{(j)}_Y \equiv 0$, $j=1,2, \cdots, n-1$; $Y=U,V$; then
\[
M^{(n)}_Y = \int^{\infty}_{-\infty} \langle \na Y(\omega^{(0)}) \, , 
\, h^{(n)} \rangle \, dt \, , \ \ Y=U,V\ ;
\]
where $h^{(n)}$ is given in (\ref{hva}) and $\om^{(j)}$ ($j=1,2, \cdots, 
n-1$) solves the linear equation (\ref{dif}) under the boundary 
condition (\ref{difbc}).}}

\subsection{Numerical Evaluations of the Second and Third Order Melnikov 
Functions}

\subsubsection{Numerical Evaluations of the Second Order Melnikov 
Functions}

In the formulae of the second order Melnikov functions
[(\ref{hmelf}), for $n=2$], since the functions $U$ and $V$ only depend on
$\omega_j \, (1 \leq j \leq 4)$, only the following four
components of $h^{(2)}$ enter the evaluation:
\begin{eqnarray}
  h^{(2)}_1 &=& -A_2 \omega^{(1)}_p \omega^{(1)}_2
       +  A_0 [ \omega^{(0)}_p \omega^{(1)}_0 
       + \omega^{(1)}_p \omega^{(0)}_0] \, , \label{exh1} \\[1ex]
    h^{(2)}_2 &=&  A_1 \omega^{(1)}_p \omega^{(1)}_1
       -  A_2 \omega^{(1)}_p \omega^{(1)}_3 \, , \label{exh2}\\[1ex]
    h^{(2)}_3 &=&  A_2 \omega^{(1)}_p \omega^{(1)}_2
       -  A_1 \omega^{(1)}_p \omega^{(1)}_4 \, , \label{exh3}\\[1ex]
   h^{(2)}_4 &=&   A_2 \omega^{(1)}_p \omega^{(1)}_3 
        - A_0 [ \omega^{(0)}_p \omega^{(1)}_5 
       + \omega^{(1)}_p \omega^{(0)}_5] \, . \label{exh4} 
\end{eqnarray}
Thus we need to know the functions $\omega^{(1)}_p$,
$\omega^{(1)}_j$ $(0 \leq j \leq 5)$.  Specifically these
functions satisfy a set of decoupled systems of linear
equations.  The five functions  $\omega^{(1)}_p$,
$\omega^{(1)}_{\ell}$ $(1 \leq \ell \leq 4)$ satisfy a
self-contained system of linear equations with variable
coefficients:
\begin{eqnarray}
  \label{sseq1}
  \left(
    \begin{array}{l}
       \omega^{(1)}_1 \\[1ex]
       \omega^{(1)}_2 \\[1ex]
        \omega^{(1)}_3 \\[1ex]
       \omega^{(1)}_4 \\[1ex]
       \omega^{(1)}_p \\[1ex]
    \end{array} \right)^{\bullet} = \AA \left(
    \begin{array}{l}
       \omega^{(1)}_1 \\[1ex]
       \omega^{(1)}_2 \\[1ex]
        \omega^{(1)}_3 \\[1ex]
       \omega^{(1)}_4 \\[1ex]
       \omega^{(1)}_p \\[1ex]
    \end{array} \right) + \BB \, ,
\end{eqnarray}
where
\begin{eqnarray*}
  \AA &=& \left(
    \begin{array}{ccccc}
      0 & -A_2 \omega^{(0)}_p & 0 & 0 &  -A_2 \omega^{(0)}_2 \\[1ex]
       A_1 \omega^{(0)}_p & 0 &  -A_2 \omega^{(0)}_p & 0 &
        A_1 \omega^{(0)}_1  - A_2 \omega^{(0)}_3 \\[1ex]
        0 &  A_2 \omega^{(0)}_p & 0 &  -A_1 \omega^{(0)}_p &
           A_2 \omega^{(0)}_2  - A_1 \omega^{(0)}_4 \\[1ex]
          0 & 0 &  A_2 \omega^{(0)}_p & 0 &  A_2 \omega^{(0)}_3
          \\
 -A_{1,2} \omega^{(0)}_2 &  -A_{1,2} \omega^{(0)}_1 &
  A_{1,2} \omega^{(0)}_4 &  A_{1,2} \omega^{(0)}_3 & 0
    \end{array} \right) \, , \\[1ex]
\BB &=& \left(
  \begin{array}{c}
    A_0  \omega^{(0)}_p  \omega^{(0)}_0 \\[1ex]
    0 \\[1ex]
    0 \\[1ex]
 -A_0 \omega^{(0)}_p  \omega^{(0)}_5 \\[1ex]
 -  A_{4,5} \omega^{(0)}_4 \omega^{(0)}_5 
 -  A_{0,1} \omega^{(0)}_0 \omega^{(0)}_1 
  \end{array} \right) \, .
\end{eqnarray*}
The functions $\omega^{(1)}_0$ and $ \omega^{(1)}_5 $ satisfy the 
linear equations:
\begin{eqnarray}
  \label{sseq2}
  \dot{\omega}^{(1)}_0 &=& \omega^{(0)}_p [
  A_{-1} \omega^{(1)}_{-1} - A_1 \omega^{(1)}_1 ]
  -A_1 \omega^{(0)}_1  \omega^{(1)}_p \, ,\\[1ex]
  \label{sseq3}
  \dot{\omega}^{(1)}_5 &=& \omega^{(0)}_p [
  A_{4} \omega^{(1)}_{4} - A_6 \omega^{(1)}_6 ]
  + A_4 \omega^{(0)}_4  \omega^{(1)}_p \, .
\end{eqnarray}
The two functions $\omega^{(1)}_{-1}$ and  $\omega^{(1)}_6$ on
the right hand sides of (\ref{sseq2}) and (\ref{sseq3}) are
obtained from solving the following two self-contained systems of 
linear equations:
\begin{eqnarray}
  \label{sseq4}
\frac{d}{d \zeta} \left[
  \begin{array}{c}
\omega^{(1)}_{-1} \\[1ex]
\omega^{(1)}_{-2} \\[1ex]
\omega^{(1)}_{-3} \\[1ex]
\omega^{(1)}_{-4} \\[1ex]
\end{array} \right] = \C \left[
  \begin{array}{c}
\omega^{(1)}_{-1} \\[1ex]
\omega^{(1)}_{-2} \\[1ex]
\omega^{(1)}_{-3} \\[1ex]
\omega^{(1)}_{-4} \\[1ex]
  \end{array} \right] + \DD \, , \\[1ex]
  \label{sseq5}
\frac{d}{d \zeta} \left[
  \begin{array}{c}
\omega^{(1)}_{6} \\[1ex]
\omega^{(1)}_{7} \\[1ex]
\omega^{(1)}_{8} \\[1ex]
\omega^{(1)}_{9} \\[1ex]
  \end{array} \right] = \EE \left[
  \begin{array}{c}
\omega^{(1)}_{6} \\[1ex]
\omega^{(1)}_{7} \\[1ex]
\omega^{(1)}_{8} \\[1ex]
\omega^{(1)}_{9} \\[1ex]
  \end{array} \right] + \FF \, , 
\end{eqnarray}
where $\frac{d}{d \zeta} = \frac{1}{\omega^{(0)}_p} \frac{d}{dt}$, 
$\zeta = \frac{1}{\k}
\ \mbox{ln}\ \cosh ( \k \Gamma t + \tau_0)$,
\begin{eqnarray*}
  \C &=& \left[
    \begin{array}{cccc}
      0 & A_{-2} & 0 & 0 \\[1ex]
      -A_{-1} & 0 & A_{-3} & 0 \\[1ex]
      0 & -A_{-2} & 0 & A_{-4} \\[1ex]
      0 & 0 & -A_{-3} & 0
    \end{array} \right] \\[1ex]
&=& \left[ 
    \begin{array}{cccc}
      0 & -\frac{39}{82} & 0 & 0 \\[1ex]
      \frac{23}{50} & 0 & - \frac{59}{122} & 0 \\[1ex]
      0 & \frac{39}{82} & 0 & - \frac{83}{170} \\[1ex]
      0 & 0 & \frac{59}{122} & 0 
    \end{array} \right] \, , \\[2ex]
\DD &=&
\left[
  \begin{array}{ccccc}
    -A_0 \omega^{(0)}_0 \\[1ex]
    0 \\[1ex]
    0\\[1ex]
    0
  \end{array} \right] = \left[
  \begin{array}{ccccc}
   \frac{11}{26} \omega^{(0)}_0 \\[1ex]
    0 \\[1ex]
    0\\[1ex]
    0
  \end{array} \right] \, . \\[2ex]
\EE &=& \left[
  \begin{array}{ccccc}
    0 & -A_7 & 0 & 0 \\[1ex]
    A_6 & 0 & -A_8 & 0 \\[1ex]
    0 & A_7 & 0 & -A_9 \\[1ex]
    0 & 0 & A_8 & 0
  \end{array} \right] = - \C \, , \\[2ex]
\FF &=& \left[
  \begin{array}{c}
    A_5 \omega^{(0)}_5 \\[1ex]
    0\\[1ex]
    0\\[1ex]
    0
  \end{array} \right] = \left[
  \begin{array}{c}
    - \frac{11}{26} \omega^{(0)}_5 \\[1ex]
    0\\[1ex]
    0\\[1ex]
    0
  \end{array} \right] \, .
\end{eqnarray*}
The characteristic equation for the matrix $\C$:
\begin{displaymath}
  \lambda^4 + \lambda^2 \left[ \frac{23 \times 39}{50 \times 82} +
    \frac{59}{122} \left( \frac{39}{82} + \frac{83}{170} \right) \right] +
  \frac{23 \times 39 \times 59 \times 83}{50 \times 82 \times 122 
    \times 170} =0 \, ,
\end{displaymath}
has only imaginary eigenvalues,
\begin{displaymath}
  \lambda_{1,2} \dot{=} \pm i\ 0.773 \, , \quad
   \lambda_{3,4} \dot{=} \pm i \ 0.295 \, .
\end{displaymath}
Let $v_j \ ( 1 \leq j \leq 4)$ be the corresponding eigenvectors,
$v=$ columns $\left\{ v_1,v_2,v_3,v_4 \right\}$,\break $w^-= (
\omega^{(1)}_{-1}, \omega^{(1)}_{-2}, \omega^{(1)}_{-3},
\omega^{(1)}_{-4})^T$, and $w^+= (\omega^{(1)}_6, \omega^{(1)}_7,
\omega^{(1)}_8, \omega^{(1)}_9 )^T$; then
\begin{eqnarray}
  w^-(\zeta) &=& v\int^{\zeta}_{\sign(\k) \infty} \diag \bigg \{
    e^{\lambda_1(\zeta - \xi)} \, , \,  e^{\lambda_2(\zeta -
      \xi)} \, , \,  e^{\lambda_3(\zeta - \xi)} \, , \,
     e^{\lambda_4(\zeta - \xi)} \bigg \} \nonumber \\
   & & \ \ \ \ \ v^{-1} \DD(\xi) \ d\xi \ , \label{soln1}\\
    w^+(\zeta) &=& v \int^{\zeta}_{\sign(\k) \infty} \diag \bigg \{
    e^{\lambda_1(\xi - \zeta)} \, , \,  e^{\lambda_2(
      \xi- \zeta)} \, , \,  e^{\lambda_3(\xi - \zeta)} \, , \,
     e^{\lambda_4(\xi -\zeta)} \bigg \} \nonumber \\
  & &  \ \ \ \ \ v^{-1} \FF(\xi) \ d\xi \ . \label{soln2} 
\end{eqnarray}
The solutions $w^{\pm}$ given above satisfy the boundary
condition $w^{\pm}(t=- \infty)=0$.  In fact, $w^{\pm}$ vanish at 
both positive and negative infinities, $w^{\pm}(t=\infty)=0$.

The system of coupled linear equations 
(\ref{sseq1}),  (\ref{sseq4}) and  (\ref{sseq5}) were 
integrated numerically by Thomas Witelski using a fourth order Runge-Kutta 
scheme. The
resulting numerical solutions were then used to evaluate the integrals
for the second order Melnikov function $M_U^{(2)}$ and $M_V^{(2)}$ 
using the trapezoidal rule. For
calculations with sufficiently small step-sizes, the
second order Melnikov functions were found to be identically zero to
machine precision, independent of the value of the parameters $\theta_0$
and $\tau_0$. Evaluation of the partial integrals for $-\infty<\tau\le t$
strongly suggest that the integrands of $M_U^{(2)}$ and $M_V^{(2)}$ 
are odd functions.

\subsubsection{Numerical Evaluations of the Third Order Melnikov 
Functions}

In the formulae of the third order Melnikov functions
[(\ref{hmelf}), for $n=3$], since the functions $U$ and $V$ only depend on
$\omega_j \, (1 \leq j \leq 4)$, only the following four
components of $h^{(3)}$ enter the evaluation:
\begin{eqnarray*}
h^{(3)}_1 &=& \bigg [ A_0 \omega^{(0)}_p \bigg ]\omega^{(2)}_0 
              - \bigg [ A_2 \omega^{(1)}_p \bigg ]\omega^{(2)}_2 \\
          & & +\bigg [  A_0 \omega^{(0)}_0 - A_2 \omega^{(1)}_2 \bigg ]
              \omega^{(2)}_p + \bigg [  A_0 \omega^{(1)}_0\omega^{(1)}_p 
              \bigg ] \ , \\
h^{(3)}_2 &=& \bigg [ A_1 \omega^{(1)}_p \bigg ]\omega^{(2)}_1 
              - \bigg [ A_2 \omega^{(1)}_p \bigg ]\omega^{(2)}_3 \\
          & & +\bigg [  A_1 \omega^{(1)}_1 - A_2 \omega^{(1)}_3 \bigg ]
              \omega^{(2)}_p \ , \\
h^{(3)}_3 &=& \bigg [ A_2 \omega^{(1)}_p \bigg ]\omega^{(2)}_2 
              - \bigg [ A_1 \omega^{(1)}_p \bigg ]\omega^{(2)}_4 \\
          & & +\bigg [  A_2 \omega^{(1)}_2 - A_1 \omega^{(1)}_4 \bigg ]
              \omega^{(2)}_p \ , \\
h^{(3)}_4 &=& \bigg [ A_2 \omega^{(1)}_p \bigg ]\omega^{(2)}_3 
              - \bigg [ A_5 \omega^{(0)}_p \bigg ]\omega^{(2)}_5 \\
          & & +\bigg [  A_2 \omega^{(1)}_3 - A_5 \omega^{(0)}_5 \bigg ]
              \omega^{(2)}_p - \bigg [  A_5 \omega^{(1)}_5\omega^{(1)}_p 
              \bigg ] \ .
\end{eqnarray*}
Thus we need to know the functions $\omega^{(2)}_p$,
$\omega^{(2)}_j$ $(0 \leq j \leq 5)$.  Specifically these
functions satisfy a set of decoupled systems of linear
equations.  The five functions  $\omega^{(2)}_p$,
$\omega^{(2)}_{\ell}$ $(1 \leq \ell \leq 4)$ satisfy a
self-contained system of linear equations with variable
coefficients:
\[
  \left(
    \begin{array}{l}
       \omega^{(2)}_1 \\[1ex]
       \omega^{(2)}_2 \\[1ex]
        \omega^{(2)}_3 \\[1ex]
       \omega^{(2)}_4 \\[1ex]
       \omega^{(2)}_p \\[1ex]
    \end{array} \right)^{\bullet} = \AA \left(
    \begin{array}{l}
       \omega^{(2)}_1 \\[1ex]
       \omega^{(2)}_2 \\[1ex]
        \omega^{(2)}_3 \\[1ex]
       \omega^{(2)}_4 \\[1ex]
       \omega^{(2)}_p \\[1ex]
    \end{array} \right) + \BB^{(2)} \, ,
\]
where $\AA$ is given in (\ref{sseq1}), and
\[
\BB^{(2)} = \left(
  \begin{array}{c}
     h_1^{(2)}\\[1ex]  h_2^{(2)}\\[1ex]  h_3^{(2)}\\[1ex] h_4^{(2)}\\[1ex]
  h_p^{(2)}\\[1ex]
\end{array} \right) \, ,
\]
where $h_1^{(2)}$,  $h_2^{(2)}$, $h_3^{(2)}$ and $h_4^{(2)}$ are given in 
(\ref{exh1}), (\ref{exh2}), (\ref{exh3}) and (\ref{exh4}), and 
\begin{eqnarray*}
h_p^{(2)} &=& - \bigg [  A_{-4,-3} \omega^{(1)}_{-4}\omega^{(1)}_{-3} +  
A_{-3,-2} \omega^{(1)}_{-3}\omega^{(1)}_{-2} + A_{-2,-1} 
\omega^{(1)}_{-2}\omega^{(1)}_{-1} \\ 
& & + A_{-1,0} \omega^{(1)}_{-1}\omega^{(0)}_{0} + A_{0,1} (\omega^{(0)}_{0}
\omega^{(1)}_{1} + \omega^{(1)}_{0}\omega^{(0)}_{1}) \\
& & + A_{1,2} \omega^{(1)}_{1}\omega^{(1)}_{2} +  
A_{2,3} \omega^{(1)}_{2}\omega^{(1)}_{3} + A_{3,4} 
\omega^{(1)}_{3}\omega^{(1)}_{4} \\ 
& & + A_{4,5} (\omega^{(0)}_{4}
\omega^{(1)}_{5} + \omega^{(1)}_{4}\omega^{(0)}_{5})+ A_{5,6} 
\omega^{(0)}_{5}\omega^{(1)}_{6} \\
& & + A_{6,7} \omega^{(1)}_{6}\omega^{(1)}_{7} +  
A_{7,8} \omega^{(1)}_{7}\omega^{(1)}_{8} + A_{8,9} 
\omega^{(1)}_{8}\omega^{(1)}_{9} \bigg ] \ .
\end{eqnarray*}
This expression is obtained from the fact that $\omega^{(1)}_{l} = 0$
for $l < - 5$ and $l > 10$. The two functions $\omega^{(2)}_0$ and 
$\omega^{(2)}_5 $ satisfy the linear equations:
\begin{eqnarray}
  \dot{\omega}^{(2)}_0 &=& \bigg [ A_{-1}\omega^{(0)}_p \bigg ] 
\omega^{(2)}_{-1} - \bigg [ A_{1}\omega^{(0)}_p \bigg ] 
\omega^{(2)}_{1} - \bigg [ A_{1}\omega^{(0)}_1 \bigg ] 
\omega^{(2)}_{p} + h^{(2)}_0 \ , \label{tseq2} \\[1ex]
  \dot{\omega}^{(2)}_5 &=& \bigg [ A_{4}\omega^{(0)}_p \bigg ] 
\omega^{(2)}_{4} + \bigg [ A_{4}\omega^{(0)}_4 \bigg ] 
\omega^{(2)}_{p} - \bigg [ A_{6}\omega^{(0)}_p \bigg ] 
\omega^{(2)}_{6} + h^{(2)}_5 \ , \label{tseq3}
\end{eqnarray}
where 
\begin{eqnarray*}
h^{(2)}_0 &=& A_{-1}\omega^{(1)}_p \omega^{(1)}_{-1} - 
A_{1}\omega^{(1)}_p \omega^{(1)}_{1} \ , \\
h^{(2)}_5 &=& A_{4}\omega^{(1)}_p \omega^{(1)}_{4} - 
A_{6}\omega^{(1)}_p \omega^{(1)}_{6} \ .
\end{eqnarray*}
The two functions $\omega^{(2)}_{-1}$ and  $\omega^{(2)}_6$ on
the right hand sides of (\ref{tseq2}) and (\ref{tseq3}) are
obtained from solving the following two self-contained systems of 
linear equations:
\begin{eqnarray}
\frac{d}{d \zeta} \left[
  \begin{array}{c}
\omega^{(2)}_{-1} \\[1ex]
\omega^{(2)}_{-2} \\[1ex]
\omega^{(2)}_{-3} \\[1ex]
\omega^{(2)}_{-4} \\[1ex]
\end{array} \right] = \C \left[
  \begin{array}{c}
\omega^{(2)}_{-1} \\[1ex]
\omega^{(2)}_{-2} \\[1ex]
\omega^{(2)}_{-3} \\[1ex]
\omega^{(2)}_{-4} \\[1ex]
  \end{array} \right] + \DD^{(2)} \, , \label{tseq4} \\[1ex]
\frac{d}{d \zeta} \left[
  \begin{array}{c}
\omega^{(2)}_{6} \\[1ex]
\omega^{(2)}_{7} \\[1ex]
\omega^{(2)}_{8} \\[1ex]
\omega^{(2)}_{9} \\[1ex]
  \end{array} \right] = \EE \left[
  \begin{array}{c}
\omega^{(2)}_{6} \\[1ex]
\omega^{(2)}_{7} \\[1ex]
\omega^{(2)}_{8} \\[1ex]
\omega^{(2)}_{9} \\[1ex]
  \end{array} \right] + \FF^{(2)} \, , \label{tseq5}
\end{eqnarray}
where $\zeta$, $\C$ and $\EE$ are given in (\ref{sseq4}) and (\ref{sseq5}),
\[
\DD^{(2)} = \frac{1}{\om^{(0)}_p}
\left[
  \begin{array}{c}
    h_{-1}^{(2)} \\[1ex]
    h_{-2}^{(2)} \\[1ex]
    h_{-3}^{(2)} \\[1ex]
    h_{-4}^{(2)} \\
  \end{array} \right] \ ,  \ \ 
\FF^{(2)} = \frac{1}{\om^{(0)}_p}
\left[
  \begin{array}{c}
    h_{6}^{(2)} \\[1ex]
    h_{7}^{(2)} \\[1ex]
    h_{8}^{(2)} \\[1ex]
    h_{9}^{(2)} \\
  \end{array} \right] \ ,
\]
\begin{eqnarray*}
h^{(2)}_{-1} &=& A_{-2}\omega^{(1)}_p \omega^{(1)}_{-2} - 
A_{0}(\omega^{(0)}_p \omega^{(1)}_{0} +\omega^{(1)}_p \omega^{(0)}_{0}) \ , \\
h^{(2)}_{-2} &=& A_{-3}\omega^{(1)}_p \omega^{(1)}_{-3} - 
A_{-1}\omega^{(1)}_p \omega^{(1)}_{-1}\ , \\
h^{(2)}_{-3} &=& A_{-4}\omega^{(1)}_p \omega^{(1)}_{-4} - 
A_{-2}\omega^{(1)}_p \omega^{(1)}_{-2}\ , \\
h^{(2)}_{-4} &=& A_{-5}\omega^{(0)}_p \omega^{(1)}_{-5} - 
A_{-3}\omega^{(1)}_p \omega^{(1)}_{-3}\ , \\
h^{(2)}_{6} &=& A_{5}(\omega^{(0)}_p \omega^{(1)}_{5} + \omega^{(1)}_p 
\omega^{(0)}_{5}) - A_{7}\omega^{(1)}_p \omega^{(1)}_{7}\ , \\
h^{(2)}_{7} &=& A_{6}\omega^{(1)}_p \omega^{(1)}_{6} 
- A_{8}\omega^{(1)}_p \omega^{(1)}_{8}\ , \\
h^{(2)}_{8} &=& A_{7}\omega^{(1)}_p \omega^{(1)}_{7} 
- A_{9}\omega^{(1)}_p \omega^{(1)}_{9}\ , \\
h^{(2)}_{9} &=& A_{8}\omega^{(1)}_p \omega^{(1)}_{8} 
- A_{10}\omega^{(0)}_p \omega^{(1)}_{10}\ .
\end{eqnarray*}
In these formulas, $\omega^{(1)}_{-5}$ and $\omega^{(1)}_{10}$ appeared, 
and they can be obtained from,
\[
\left \{ \begin{array}{l} \dot{\omega}^{(1)}_{-5} = - A_{-4}\omega^{(0)}_p 
\omega^{(1)}_{-4} \ , \\ \dot{\omega}^{(1)}_{10} = A_{9}\omega^{(0)}_p 
\omega^{(1)}_{9} \ . \\ \end{array}\right.
\]
Let $\hat{w}^-= (\omega^{(2)}_{-1}, \omega^{(2)}_{-2}, \omega^{(2)}_{-3},
\omega^{(2)}_{-4})^T$, and $\hat{w}^+= (\omega^{(2)}_6, \omega^{(2)}_7,
\omega^{(2)}_8, \omega^{(2)}_9 )^T$; then we have the similar representations 
as for $w^-$ and $w^+$,
\begin{eqnarray}
  \hat{w}^-(\zeta) &=& v\int^{\zeta}_{\sign(\k) \infty} \diag \bigg \{
    e^{\lambda_1(\zeta - \xi)} \, , \,  e^{\lambda_2(\zeta -
      \xi)} \, , \,  e^{\lambda_3(\zeta - \xi)} \, , \,
     e^{\lambda_4(\zeta - \xi)} \bigg \} \nonumber \\
   & & \ \ \ \ \ v^{-1} \DD^{(2)}(\xi) \ d\xi \ , \label{tsoln1}\\
    \hat{w}^+(\zeta) &=& v \int^{\zeta}_{\sign(\k) \infty} \diag \bigg \{
    e^{\lambda_1(\xi - \zeta)} \, , \,  e^{\lambda_2(
      \xi- \zeta)} \, , \,  e^{\lambda_3(\xi - \zeta)} \, , \,
     e^{\lambda_4(\xi -\zeta)} \bigg \} \nonumber \\
  & &  \ \ \ \ \ v^{-1} \FF^{(2)}(\xi) \ d\xi \ . \label{tsoln2} 
\end{eqnarray}
The solutions $\hat{w}^{\pm}$ given above satisfy the boundary
condition $\hat{w}^{\pm}(t=- \infty)=0$.  In fact, $\hat{w}^{\pm}$ vanish at 
both positive and negative infinities, $\hat{w}^{\pm}(t=\infty)=0$.

The above systems are solved under the boundary condition,
\[
\omega^{(2)}_l(t=-\infty)=0\ , \ \ \forall l \in Z\ .
\]
Moreover, direct verification shows that $\omega^{(2)}_l = 0$ for 
$l < -10$ and $l > 15$. For this case, the numerics of Thomas Witelski 
could not give a firm conclusion, and we are continuing investigating 
the numerics. This numerics is a part of a future project on numerical 
investigation of the degeneracy v.s. nondegeneracy of the hyperbolic 
foliations of the 2D Euler equation.

\eqnsection{Conclusion and Discussion}

In this Part II of our study, we see an explicit representation 
for the degeneracy of the hyperbolic structures for a Galerkin truncation. 
We use higher order Melnikov functions to study the robustness of such 
structures for the so-called dashed-line model. We find that both the 
first and the second order Melnikov functions are identically zero, 
which indicates that such structures are relatively robust. The study 
in this paper serves as a clue in searching for homoclinic structures 
for 2D Euler equation. The recent breakthrough result \cite{Li00c} of 
mine on the existence of a Lax pair for 2D Euler equation makes it 
plausible for the existence of homoclinic structures. \\

{\bf Acknowlegment:} The author is greatly indebted for the numerical 
contribution of Professor Thomas Witelski. The work of Thomas was 
carried at MIT and Duke University. The author 
was surprised in finding out that Thomas did not want to be a co-author,
rather gave all the credits to the author.

\newpage
\bibliography{dash}

\end{document}